\documentclass[final]{siamltex}
\usepackage{algorithm,algpseudocode}
\usepackage{latexsym, graphicx, epsfig, amsmath, amsfonts,amssymb,bm}
\usepackage{caption, subcaption}
\usepackage{epstopdf}
\usepackage{booktabs}
\usepackage{array}
\usepackage{color}
\newtheorem{remark}{Remark}[section]

\def\be{\begin{equation}}
\def\ee{\end{equation}}

\def\x{\mathbf{x}}
\def\y{\mathbf{y}}
\def\f{\mathbf{f}}

\def\z{\mathbf{z}}
\def\PPhi{\boldsymbol{\Phi}}
\def\pphi{\boldsymbol{\phi}}
\def\N{\mathbf{N}}
\def\W{\mathbf{W}}

\def\R{{\mathbb R}}

\title{Data Driven Governing Equations Approximation Using Deep Neural Networks}

\author{Tong Qin\and Kailiang Wu\and Dongbin
       Xiu\thanks{Department of Mathematics,
		The Ohio State University, Columbus, OH 43210, USA.
		{\tt qin.428@osu.edu, wu.3423@osu.edu, xiu.16@osu.edu.}
		Funding: This work was partially supported by AFOSR FA9550-18-1-0102.}
}

\begin{document}
\maketitle
\begin{abstract}
We present a numerical framework for approximating unknown governing equations using observation data and deep neural networks (DNN).  In particular, we propose to use residual network (ResNet) as the basic building block for equation approximation. We demonstrate that the ResNet block can be considered as a one-step method that is exact in temporal integration. We then 
present two multi-step methods, recurrent ResNet (RT-ResNet) method and recursive ReNet (RS-ResNet) method.  The
RT-ResNet is a multi-step method on uniform time steps, whereas the RS-ResNet is an adaptive 
multi-step method using variable time steps. All three methods presented here are based on integral form of the underlying
dynamical system. As a result, they do not require time derivative data for equation recovery and can cope with relatively coarsely distributed trajectory data. Several numerical examples are presented to demonstrate the performance of the methods.
\end{abstract}

% main text
\section{Introduction} \label{sec:intro}

Recently there has been a growing interest in discovering governing
equations numerically using observational data.
Earlier efforts include methods using symbolic regression
(\cite{bongard2007automated,schmidt2009distilling}), equation-free
modeling \cite{kevrekidis2003equation}, heterogeneous multi-scale
method (HMM) (\cite{E_HMM03}), %data series \cite{crutchfield1987equations}, 
artificial neural networks (\cite{gonzalez1998identification}),
nonlinear regression (\cite{voss1999amplitude}), 
%normal form identification (\cite{majda2009normal}), 
empirical dynamic modeling (\cite{sugihara2012detecting,ye2015equation}), 
nonlinear Laplacian spectral analysis (\cite{giannakis2012nonlinear}), 
%modeling emergent behavior (\cite{roberts2014model}),  
automated inference of dynamics
(\cite{schmidt2011automated,daniels2015automated,daniels2015efficient}),
etc.
More recent efforts start to cast the problem into a function
approximation problem, where the
unknown governing equations are treated as target functions relating the data for the state variables and
their time derivatives. % is the sought-after governing equation.
The majority of the methods employ certain sparsity-promoting
algorithms to create parsimonious models from a large set of
dictionary for all possible models, so that the
true dynamics could be recovered exactly
(\cite{tibshirani1996regression}).
Many studies have been conducted to effectively deal with noises in data
(\cite{brunton2016discovering, schaeffer2017sparse}), corruptions in data
(\cite{tran2017exact}), 
partial differential equations \cite{rudy2017data,
  schaeffer2017learning}, etc.
Methods have also been developed in conjunction with
model selection approach
(\cite{Mangan20170009}), Koopman theory (\cite{brunton2017chaos}), and
Gaussian process regression (\cite{raissi2017machine}), to name a few.
%Three sampling strategies 
%were developed in \cite{schaeffer2017extracting} 
%for recovering quadratic high-dimensional differential equations from under-sampled data. 
A more recent work resorts to the more traditional means of
approximation by using orthogonal polynomials
(\cite{WuXiu_JCPEQ18}). The approach seeks accurate numerical approximation to the underlying
governing equations, instead of their exact recovery. By doing so,
many existing results in polynomial approximation theory can be
applied, particularly those on sampling strategies. It was shown in
\cite{WuXiu_JCPEQ18} that data from a large number of short bursts of
trajectories are more effective for equation recovery than those from
a single long trajectory.

On the other hand, artificial neural network (ANN), and particularly deep neural network (DNN), has seen tremendous
successes in many different disciplines. The number of publications is
too large to mention. Here we cite only a few relatively more recent
review/summary type publications \cite{montufar2014number, 
bianchini2014complexity, eldan2016power, poggio2017, DuSwamy2014,
GoodfellowBC-2016, Schmidhuber2015}.  
Efforts have been devoted to the use of ANN for
various aspects of scientific computing, including construction of
reduced order model (\cite{HesthavenU_JCP18}), aiding solution of
conservation laws (\cite{RayHeasthaven_JCP18}), multiscale problems
(\cite{ChanE_JCP18, Wang_2018}),
solving and learning systems involving ODEs and PDEs
(\cite{MardtPWN_Nature18, Chen_2018, long2017pde, KhooLuYing_2018}),
uncertainty quantification (\cite{TripathyB_JCP18, Zabaras_2018}), etc.

The focus of this paper is on the approximation/learning of dynamical
systems using deep neural networks (DNN). The topic has been
explored in a series of recent articles, in the context of ODEs
\cite{raissi2018multistep,rudy2018deep}) and PDEs
(\cite{raissi2018hidden,raissi2018deep,long2017pde}).
The new contributions of this paper include the following. First,
we introduce new constructions of
deep neural network (DNN), specifically suited for learning dynamical
systems. In particular, our new network structures employ residual
network (ResNet), which was first proposed in \cite{he2016deep} for
image analysis and has become very popular due to its effectiveness.
In our construction, we employ a ResNet block, which consists of
multiple fully connected hidden layers, as the fundamental building block of our
DNN structures. We show that the ResNet block can be considered as
a one-step numerical integrator in time. This integrator is ``exact''
in time, i.e., no temporal error, in the sense that the only error stems from the neural
network approximation of the evolution operators defining the governing equation. This is different
from a few existing work where ResNet is viewed as the Euler forward
scheme (\cite{ChangEtAl2018}).
Secondly, we introduce two variations of the ResNet structure to serve as
multi-step learning of the underlying governing equations. The first one employs
recurrent use of the ResNet block. This is inspired by the well known
recurrent neural network (RNN), whose connection with dynamical
systems has long been recognized, cf. \cite{GoodfellowBC-2016}. Our
recurrent network, termed RT-ResNet hereafter, is different in the sense that the recurrence is
enforced blockwise on the ResNet block, which by itself is a DNN. (Note that in
the traditional RNN, the recurrence is enforced on the hidden layers.)
We show that
the RT-ResNet is a multi-step integrator that is exact in time, with
the only error stemming from the ResNet approximation of the evolution
operator of the 
underlying equation. The other variation of the ResNet approximator
employs recursive use of the ResNet block, termed RS-ResNet. Again,
the recursion is enforced blockwise on the ResNet block (which is a
DNN). We show
that the RS-ResNet is also an exact multi-step integrator. The
difference between RT-ResNet and RS-ResNet is that the former is
equivalent to a multi-step integrator using an uniform time step, whereas 
the latter is an ``adaptive'' method with variable time steps
depending on the particular
problem and data.
Thirdly, the derivations in this paper utilize integral form of
the underlying dynamical system. By doing so, the proposed methods do not
require knowledge or data of the time derivatives of the equation
states. This is different from most of the 
existing studies
(cf. \cite{bongard2007automated,brunton2016discovering,schaeffer2017sparse,WuXiu_JCPEQ18}),
which deal with the equations directly and thus require time
derivative data. Acquiring time derivatives introduces an additional
source for noises
and errors, particularly when one has to conduct numerical
differentiation of noisy trajectory data.
%
%
%is applied to ResNet blocks, rather than individual hidden layers.  
%This is motivated by the key observation that a long-time flow map
%can be  expressed the composition of many short-time flow maps,  
%which are near-identity maps and can be approximated efficiently by ResNet.
%
Consequently, the proposed three new DNN structures, the one-step ResNet
and multi-step RT-ResNet and RS-ResNet, are capable of approximating
unknown dynamical systems using only state variable data, which could
be relatively coarsely distributed in time.
In this case, 
most of the existing methods become less effective, as accurate
extration of time derivatives is difficult.

This paper is organized as follows. After the basic problem setup in
Section \ref{sec:setup}, we present the main methods in Section
\ref{sec:method} and some theoretical properties in Section \ref{sec:theory}. We then present, in Section
\ref{sec:examples}, a set of numerical examples, covering both linear
and nonlinear differential equations, to
demonstrate the effectiveness of the proposed algorithms.

\section{Setup} \label{sec:setup}

Let us consider an autonomous system
\be \label{govern}
\frac{d\x}{dt} = \f(\x), \qquad \x(t_0) = \x_0,
\ee
where $\x\in\R^n$ are the state variables. 
Let $\PPhi: \R^n\to \R^n$ be the flow map. The solution can be
written as
\be
\x(t; \x_0, t_0) = \PPhi_{t-t_0}(\x_0).
\ee
Note that for autonomous systems the time variable $t$ can be
arbitrarily shifted and only the time difference, or time lag, $t-t_0$ is
relevant. Hereafter we will omit $t$ in the exposition, unless
confusion arises.

In this paper, we assume  the form of the governing equations
$\f:\R^n\to\R^n$ is unknown. Our goal is to
create an accurate model for the governing equation using data
of the solution trajectories. In particular, we assume data are
collected in the form of pairs, each of which corresponds to the solution states
along one trajectory at two different time instances. That is, we
consider the set
\be \label{set}
\mathcal{S} = \{ (\z_j^{(1)}, \z_j^{(2)}): j=1,\dots, J\},
\ee
where $J$ is the total number of data pairs, and for each pair $j=1,\dots, J$,
\be
\z_j^{(1)}  = \x_j + \epsilon^{(1)}_j, \qquad \z_j^{(2)} = \PPhi_{\Delta_j}(\x_j) +\epsilon^{(2)}_j.
\ee
Here the terms $\epsilon^{(1)}_j$ and $\epsilon^{(2)}_j$ stand for the
potential noises in the data, and $\Delta_j$ is the time lag
between the two states. For notational convenience, we assume
$\Delta_j = \Delta$ to be a constant for all $j$ throughout this
paper. Consequently, the data set becomes input-output measurements of the $\Delta$-lag flow
map,
\be \label{Phi}
\x\to \PPhi_\Delta (\x).
\ee

\section{Deep Neural Network Approximation} \label{sec:method}

The core building block of our methods is a standard fully connected
feedforward neural network (FNN) with $M\geq 3$ layers, of which $(M-2)$ are hidden layers.
%
%For notational convenience and without loss of generality, we assume
%each hidden layer has the same number $m>1$ neurons, whereas both the
%input and output layers have $n$ units.
It has been established that fully connected
FNN can approximate arbitrarily well a large class of input-output maps,
i.e., they are universal approximators,  cf.~\cite{pinkus1999,
  barron1993universal, hornik1991approximation}.
Since the right-hand-side $\f$ of
\eqref{govern} is our approximation goal, we will consider $\R^n\to\R^n$ map. Let $n_j$,
$j=1, \dots, M$, be the number of neurons in each layer, we then have $n_1=n_M=n$.

Let $\N:\R^n\to\R^n$ be the operator of this network.
For any input $\y^{in}\in \R^n$, the
output of the network is 
\be \label{NN}
\y^{out} = \N (\y^{in}; \Theta),
\ee
where $\Theta$ is the parameter set including all the parameters in
the network. The operator $\N$ is a composition of the following operators
\be
\N(\cdot; \Theta) =  (\sigma_M\circ\W_{M-1})\circ \cdots\circ
(\sigma_2\circ\W_1),
\ee
%Here $\Theta=\{\W_1\dots, \W_N\}$ stands for the set of all the
%parameters in the network, 
where $\W_j$ is a matrix for the weight parameters connecting the neurons from $j$-th layer
to $(j+1)$-th layer, after using the standard approach of augmenting the
biases into the weights. The activation function $\sigma_j:\R\to \R$ is
applied component-wise to the $j$-th layer. There exist many choices for the
activation functions, e.g., sigmoid functions, ReLU (rectified linear
unit), etc. In this paper we use a sigmoid function, in particular, the $\sigma_i(x)=\tanh(x)$ function, in all layers, except
at the output layer $\sigma_M(x) = x$. This is one of the common choices
for DNN.

Using the data set \eqref{set}, we can directly train \eqref{NN} to
approximate the $\Delta$-lag flow map \eqref{Phi}. This can be done by
applying \eqref{NN} with $\y^{in}_j = \z^{(1)}_j$ to obtain
$\y^{out}_j$ for each $j=1,\dots, J$, and then minimizing following mean squared loss function
\be \label{loss}
L(\Theta) = \frac{1}{J}\sum_{j=1}^J \left\|\y^{out}_j - \z^{(2)}_j\right\|^2,
\ee
where $\| \cdot \|$ denotes vector 2-norm hereafter.
With a slight abuse of notation, hereafter we will write $\y^{in} =
\z^{(1)}$ to stand for $\y^{in}_j=\z^{(1)}_j$ for all sample
data $j=1,\dots,
J$, unless confusion arises otherwise.

%Methods along this line of approached have been explored in recent literature.

\subsection{One-step ResNet Approximation}
\label{sec:one-step}

We now present the idea of using residual neural network (ResNet) as a
one-step approximation method.
The idea of ResNet is to explicitly introduce the identity operator in
the network and force the network to effectively approximate the
``residue'' of the input-output map. Although mathematically
equivalent, this simple transformation has been shown to be highly
advantageous in practice and become increasingly popular, after its
formal introduction in \cite{he2016deep}.

\begin{figure}[ht] 
    \centering
	\includegraphics[width=10cm]{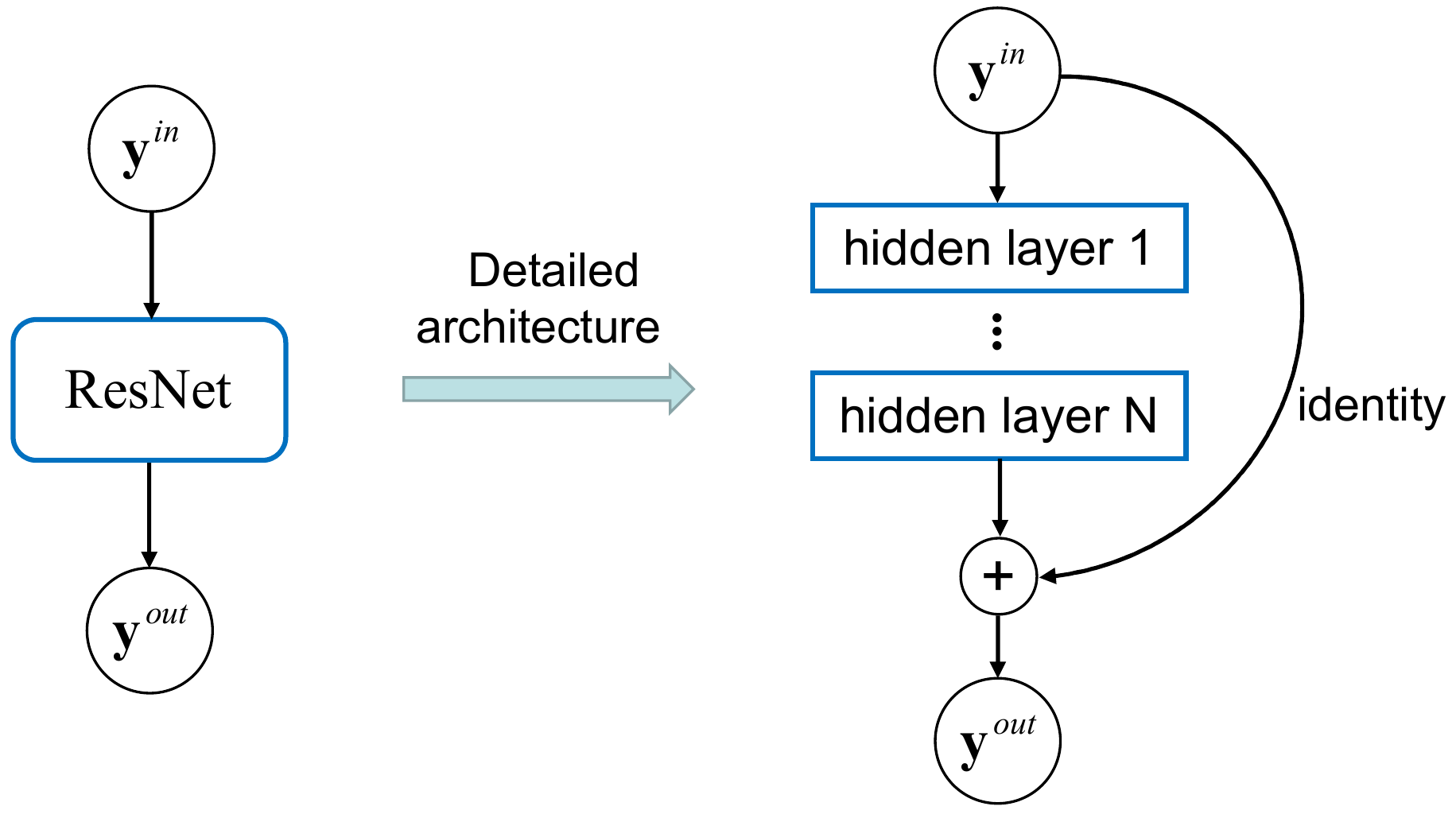} 
  \caption{Schematic of the ResNet structure for one-step approximation.}
\label{fig:ResNet}
\end{figure}

The structure of the ResNet is illustrated in Figure
\ref{fig:ResNet}. The ResNet block consists of $N$ fully connected hidden layers and
an identity operator to re-introduce the input $\y^{in}$ back into the
output of the hidden layers. This effectively produces the following mapping
\be \label{Res}
\y^{out} = \y^{in} + \N (\y^{in}; \Theta),  \quad \y^{in} = \z^{(1)},
\ee
where $\Theta $ are the weight and bias parameters in the network. The
parameters are determined by minimizing the same loss function
\eqref{loss}. This effectively accomplishes the training of the operator $\N (\cdot; \Theta)$.
%%%%%%%%%%%%%%%%%%

The connection between dynamical systems and ResNet has been
recognized. In fact, ResNet has been viewed as the Euler forward time
integrator (\cite{ChangEtAl2018}).
To further examine its property, let us
consider the exact $\Delta$-lag flow map,
\be
\begin{split}
\x(\Delta ) &= \PPhi_\Delta (\x(0)) \\
&= \x(0) + \int_0^{\Delta} \f(\x(t)) dt \\
&= \x(0) + \Delta  \cdot \f(\x(\tau))  \\
& = \x(0) + \Delta \cdot \f(\PPhi_{\tau}(\x)), \quad 0\leq \tau\leq \Delta.
\end{split}
\ee
This is a trivial derivation using the mean value theorem. For
notational convenience, we now define ``effective increment''.
\begin{definition} \label{def}
For a given autonomous system \eqref{govern}, given an initial state
$\x$ and an increment $\Delta
\geq 0$, then its effective increment of size $\Delta$ is defined as
\be \label{pphi}
\pphi_\Delta(\x; \f) =  \Delta \cdot \f(\PPhi_{\tau}(\x)),
\ee
for some $0\leq \tau \leq \Delta$ such that
\be \label{phi}
\x(\Delta ) = \x + \pphi_\Delta(\x; \f).
\ee
\end{definition}
Note that the effective increment $\pphi_\Delta$ depends only on
its initial state $\x$, once the governing equation $\f$ and the
increment $\Delta$ are fixed.

Upon comparing the exact state \eqref{phi} and the one-step ResNet
method \eqref{Res}, it is thus easy to see that a successfully trained
network operator $\N$ is an approximation to the effective increment $\pphi_\Delta$, i.e.,
\be \label{N_error}
\N(\x; \Theta) \approx \pphi_\Delta(\x; \f).
\ee

Since the effective increment completely determines the true solution
states on a $\Delta$ interval, we can then use the ResNet
operator to approximate the solution trajectory. That is, starting
with a given initial state $\y^{(0)} = \x$, we can time march the state
\be \label{resnet_eq}
\y^{(k+1)} = \y^{(k)} + \N(\y^{(k)}; \Theta), \qquad k=0,\dots.
\ee
This discrete dynamical system serves as our approximation to the true
dynamical system \eqref{govern}. It gives us an approximation to the true states on a uniform
time grids with stepsize $\Delta$. 

\begin{remark}
Even though the approximate system \eqref{resnet_eq} resembles the
well known Euler forward time stepping scheme, it is not a first-order
method in time. In fact, upon comparing \eqref{resnet_eq} and the true
state \eqref{phi}, it is easy to see that \eqref{resnet_eq} is
``exact'' in term of temporal integration. The
only source of error in the system \eqref{resnet_eq} is the
approximation error of the effective increment in \eqref{N_error}. The
size of this error is determined by the quality of the data and the
network training algorithm.
\end{remark}

\begin{remark}
The derivation here is based on \eqref{phi}, which is from
the integral form of the governing equation. As a result, training of
the ResNet method does not require data on the time derivatives of the
true states. Moreover, $\Delta$ does not need to be exceedingly
small (to enable accurate numerical differentiation in time). This makes the ResNet method suitable for problems with
relatively coarsely distributed data.
\end{remark}

\subsection{Multi-step Recurrent ResNet (RT-ResNet) Approximation}
\label{sec:multi-step_recurrent}

We now combine the idea of recurrent neural network (RNN) and
the ResNet method from the previous section. The distinct feature of our construction is that the
recurrence is applied to the entire ResNet block, rather than to the
individual hidden layers, as is done for the standard RNNs. (For an
overview of RNN, interested readers are referred to
\cite{GoodfellowBC-2016}, Ch. 10.)

\begin{figure}[ht] 
    \centering
	\includegraphics[width=12cm]{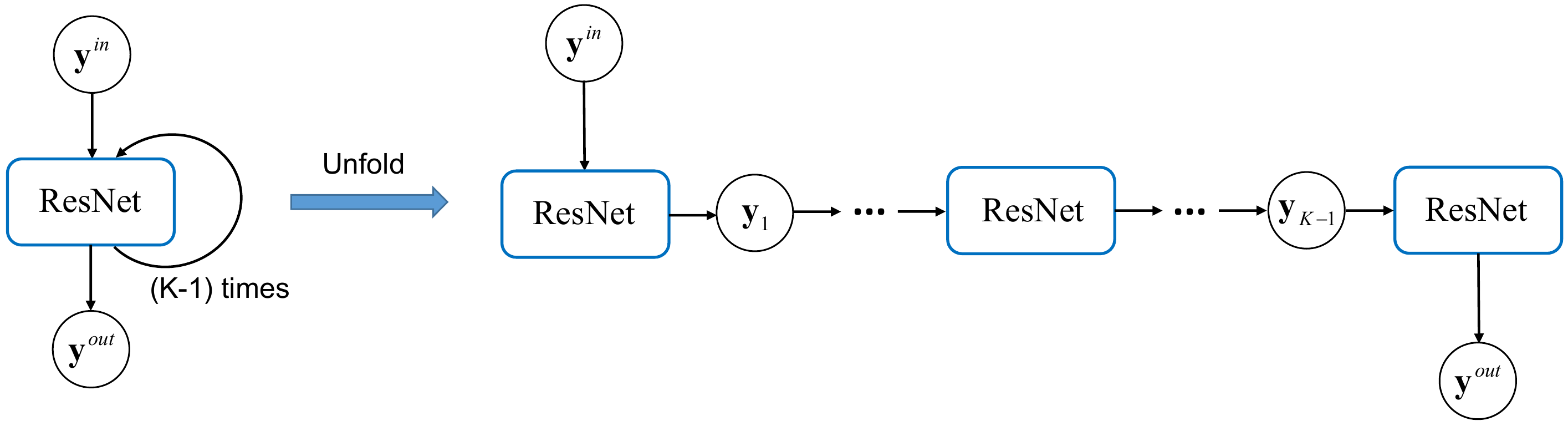} 
  \caption{Schematic of the recurrent ResNet (RT-ResNet) structure for
    multi-step approximation ($K\geq 1$).}
\label{fig:RT-ResNet}
\end{figure}

The structure of the resulting Recurrent ResNet (RT-ResNet) is shown
in Figure \ref{fig:RT-ResNet}. The ResNet block, as presented
in Figure \ref{fig:ResNet}, is ``repeated'' ($K-1$)
times, for an integer $K\geq 1$, before producing the output
$\y^{out}$. This makes the occurrence of the ResNet block a total of $K$
times.
The unfolded structure is shown on the right of
Figure \ref{fig:RT-ResNet}. The
RT-ResNet then produces the following scheme, for $K\geq 1$,
\be \label{RRes}
\left\{
\begin{split}
\y_0 &= \z^{(1)}, \\
\y_{k+1} &= \y_k + \N (\y_k; \Theta),  \quad k=0,\dots, K-1,\\
\y^{out} &= \y_K.
\end{split}
\right.
\ee
The network is then trained by using the data set \eqref{set} and
minimizing the same loss function \eqref{loss}.
%Once
%again, we remark that the recurrence is enforced on the ResNet block,
%which by itself is a deep neural network.
For $K=1$, this reduces to the one-step ResNet method \eqref{Res}.
%The
%structure of the RT-ResNet is illustrated in Figure
%\ref{fig:RT-ResNet}, using $K=3$ as an example. Note that the
%recurrence is enforced blockwise on the ResNet, which by itself is a
%DNN. The same ResNet is then recurrently used $K$ times.

%Once the network training is completed, we obtain a discrete dynamical
%system
%\be
%\y^{(k+1)} = \y^{(k)} + \N(\y^{(k)}; \Theta), \qquad k=0,\dots,
%\ee
%which gives us solution states on a discrete time with the interval
%$\delta = \Delta/K$.

To examine the properties of the RT-ResNet, let us consider a unform
discretization of the time lag $\Delta$. That is, let 
$\delta = \Delta/K$, and consider, $t_k =
k\delta$, $k=0,\dots,K$.
The exact solution state $\x$ satisfies the following relation
\be 
\left\{
\begin{split}
\x(t_0) &= \x(0), \\
\x(t_{k+1}) &= \x(t_k) + \pphi_\delta (\x(t_k); \f),  \quad k=0,\dots, K-1,\\
\x(\Delta) &= \x(t_K),
\end{split}
\right.
\ee
where $\pphi_\delta(\x;\f)$ is the effective increment
defined of size $\delta$, as defined in Definition \ref{def}.

Upon comparing this with the RT-ResNet scheme \eqref{RRes}, it is easy
to see that training the RT-ResNet is equivalent to finding the
operator $\N$ to approximate the
$\delta$-effective increment,
\be
\label{eq:NN_delta}
\N(\x;\Theta) \approx \pphi_\delta(\x; \f).
\ee
Similar to the one-step ResNet method, the multi-step RT-ResNet is
also exact in time, as it contains no temporal discretization
error. The only error stems from the approximation of the
$\delta$-effective increment.

Once the RT-ResNet is successfully trained, it gives us a discrete
dynamical system \eqref{RRes} that can be further marched in time
using any initial state. This is an approximation to the true
dynamical system on uniformly distributed time instances with an
interval $\delta=\Delta/K$. Therefore, even though the training data
are given over $\Delta$ time interval, the RT-ResNet system can
produce solution states on finer time grids with a step size
$\delta \leq \Delta$ ($K\geq 1$).

\subsection{Multi-step Recursive ResNet Approximation}
\label{sec:multi-step_recursive}

\begin{figure}[ht] 
    \centering
	\includegraphics[width=12cm]{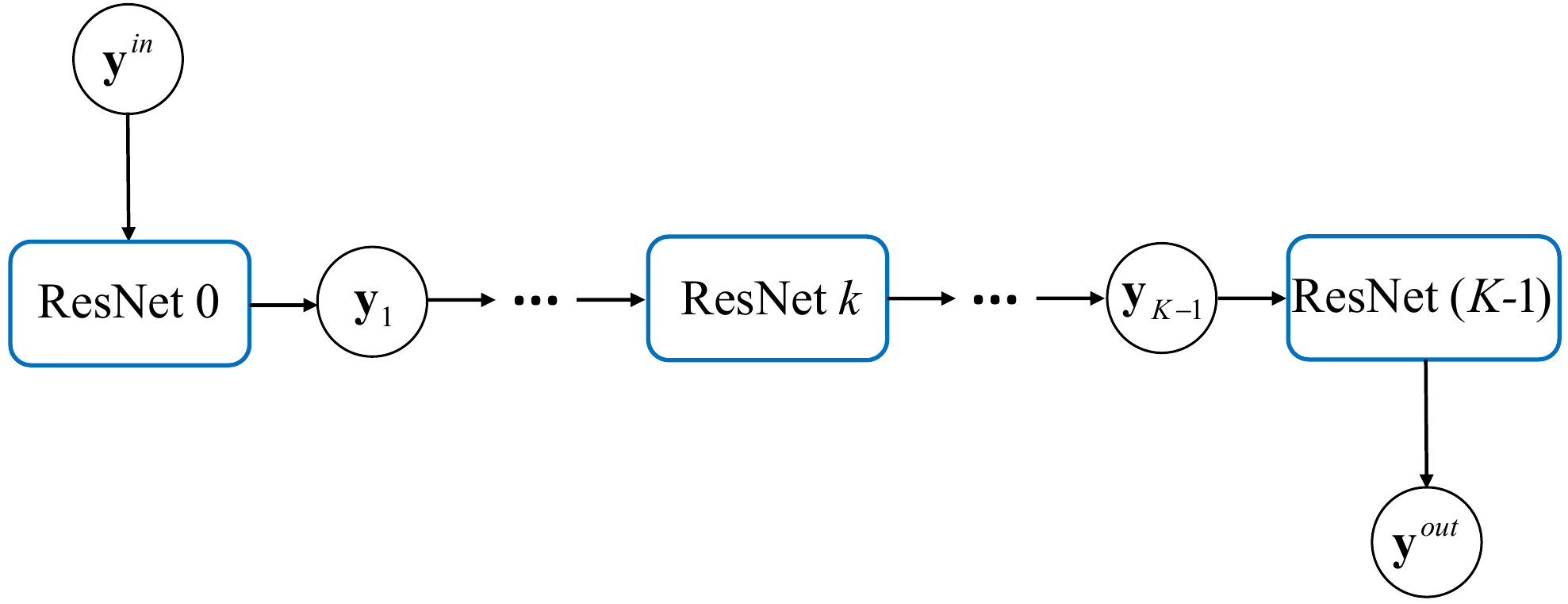} 
  \caption{Schematic of the recursive ResNet (RS-ResNet) structure for
    multi-step approximation ($K\geq 1$).}
\label{fig:RS-ResNet}
\end{figure}

We now present another multi-step approximation method based on the
ResNet block in Figure \ref{fig:ResNet}.
The structure of the network is shown in Figure
\ref{fig:RS-ResNet}. From the input $\y^{in}$, ResNet blocks are
recursively used a total of $K\geq 1$ times, before producing the
output $\y^{out}$. The network, referred to as recursive ResNet
(RS-ResNet) hereafter, thus produces the following scheme, for any
$K\geq 1$,
\be \label{Ress}
\left\{
\begin{split}
\y_0 &= \z^{(1)}, \\
\y_{k+1} &= \y_k + \N (\y_k; \Theta_k),  \quad k=0,\dots, K-1,\\
\y^{out} &= \y_K.
\end{split}
\right.
\ee
Compared to the recurrent RT-ResNet method \eqref{RRes}  from the
previous section, the major difference in RS-ResNet is that each
ResNet block inside the network has its own parameter sets $\Theta_k$ and thus are
different from each other.  Since each ResNet is a DNN by itself, the
RS-ResNet can be a very deep network when $K>1$. When $K=1$, it also
reduces back to the one-step ResNet network.

Let $0=t_0 < t_1 <\cdots < t_K = \Delta$ be an arbitrarily distributed
time instances in $[0,\Delta]$ and $\delta_k = t_{k+1} - t_k$,
$k=0,\dots, K-1$, be the (non-uniform) increments. It is then
straightforward to see that the exact state satisfies
\be 
\left\{
\begin{split}
\x(t_0) &= \x(0), \\
\x(t_{k+1}) &= \x(t_k) + \pphi_{\delta_k} (\x(t_k); \f),  \quad k=0,\dots, K-1,\\
\x(\Delta) &= \x(t_K),
\end{split}
\right.
\ee
where $\pphi_{\delta_k}(\x; \f)$ is the $\delta_k$ effective increment
defined in Definition \ref{def}.

Upon comparing with the RS-ResNet scheme \eqref{Ress}, one can see
that the training of the RS-ResNet produces the following
approximation
\be
\label{eq:RS-NN}
\N(\x;\Theta_k) \approx \pphi_{\delta_k}(\x; \f), \qquad k=0,\dots, K-1.
\ee
That is, each ResNet operator $\N(\x;\Theta_k)$ is an approximation of
an effective increment of size $\delta_k$, for $k=0,\dots, K-1$,
under the condition $\sum_{k=0}^{K-1} \delta_k = \Delta$. 
Training the
network using the data \eqref{set} and loss function \eqref{loss}
will determine the parameter sets $\Theta_k$, and subsequently the
effective increments with size $\delta_k$, for $k=0,\dots, K-1$,
From this perspective, one may view RS-ResNet as an ``adaptive''
method, as it adjusts its parameter sets to approximate $K$ smaller
effective increments whose increments are determined by the data.
Since RS-ResNet is a very deep network with a large number of
parameters, it is, in principle, capable of producing more accurate
results than ResNet and RT-ResNet, assuming cautions have been exercised to
prevent overfitting.

A successfully trained RS-ResNet also gives us a discrete dynamical
system that approximates the true governing equation
\eqref{govern}. Due to its ``adaptive'' nature, the intermediate time
intervals $\delta_k$ are variables and not known
explicitly. Therefore, the discrete RS-ResNet needs to be applied $K$
times to produce the solution states over the time interval $\Delta$,
which is the same interval given by the training data. This
is different from the RT-ResNet, which can produce solutions over a
smaller and uniform time interval $\delta =\Delta/K$.

\section{Theoretical Properties}\label{sec:theory}

In this section we present a few straightforward analysis
to demonstrate certain theoretical aspects of the proposed DNN for
equation approximation.

\subsection{Continuity of Flow Map}

Under certain conditions on ${\bf f}$, one can show that the flow map of the dynamical system 
\eqref{govern} is locally Lipschitz continuous.

\begin{lemma}\label{lem1}
Assume ${\bf f}$ is Lipschitz continuous with Lipschitz constant $L$ on a set $D \subseteq \mathbb R^n$. 
%$$
%\big\| {\bf f} ({\bf x}_1 ) - {\bf f} (  {\bf x}_2 ) \big\| \le L \| {\bf x}_1  -   {\bf x}_2  \|, \quad \forall {\bf x}_1,  {\bf x}_2 \in D.
%$$
For any $\tau>0$, define
$$
D_\tau := \Big\{ \x_0 \in D :  \PPhi_{t}( \x_0) \in D,~\forall t \in [0,  \tau]  \Big\}.
$$
Then, for any $t \in [0 , \tau]$, the flow map $\PPhi_{t}$ is Lipschitz continuous on $D_\tau$. Specifically, 
for any $\x_0, \widetilde \x_0 \in D_\tau$, 
\begin{equation}\label{eq:WKLest1}
\| \PPhi_{t}(\x_0) - \PPhi_{t}(\widetilde \x_0)  \| \le {\rm e}^{L  t } \| \x_0 - \widetilde \x_0 \|, \quad \forall t \in [0 ,   \tau].
\end{equation}
\end{lemma}

\begin{proof}
  The proof directly follows from the classical result on the continuity of the dynamical system \eqref{govern} with respect to initial data; see \cite[p. 109]{stuart1998dynamical}. 
\end{proof}

The above continuity ensures that the flow map can be approximated by neural networks to any desired degree of accuracy by increasing the number of
hidden layers and neurons; see, for example, \cite{leshno1993multilayer,pinkus1999}. 
The Lipschitz continuity will also play an important role in the error analysis in Theorem \ref{thm:err}.

\subsection{Compositions of Flow Maps}

It was shown in \cite{bartlett2018representing} that any smooth
bi-Lipschitz function
can be represented as compositions of functions, each of which is
near-identity in  
Lipschitz semi-norm.  
For the flow map of the autonomous system \eqref{govern}, 
we can prove a stronger result by using the following property
\begin{equation}\label{eq:semigroup_property}
\PPhi_{t_1} \circ \PPhi_{t_2} = \PPhi_{t_1+t_2}, \quad \forall t_1,t_2. 
\end{equation}

\begin{theorem}
	For any positive integer $K\geq 1$, the flow map $\PPhi_{\Delta} $ can be expressed as a $K$-fold composition of $\PPhi_{\delta} $, namely, 
	\begin{equation}\label{WKL1}
	\PPhi_{\Delta} = \underbrace{ \PPhi_{\delta} \circ \cdots \circ  \PPhi_{\delta} }_{{K{\rm-fold}}},
	\end{equation}
	where $\delta = \Delta/K$, and $\PPhi_{\delta} $ satisfies 
	\begin{equation}\label{WKL2}
	\|\PPhi_{\delta} (\x_0) - \x_0\| \le \frac{\Delta}{K} \sup_{t \in [0, \delta]} \| {\bf f} ( \PPhi_t (\x_0) ) \|, \quad
	\forall \x_0.
	\end{equation}
	Suppose that ${\bf f}$ is bounded on  $D \subseteq \mathbb R^n$, then 
\begin{equation}\label{WKL3}
\|\PPhi_{\delta}  - {\bf I}\|_{L^\infty(D_\delta)} \le 
 \frac{\Delta}{K}  \| {\bf f}  \|_{L^\infty (D) } = {\mathcal O} \bigg(\frac{\Delta}{K}\bigg),
\end{equation}	
where ${\bf I}: {\mathbb R}^n \to {\mathbb R}^n$ is the identity map, 
and $\| \cdot \|_{L^\infty} := {\rm ess} \sup \| \cdot \|$.
\end{theorem}

\begin{proof}
	The representation \eqref{WKL1} is a direct consequence of  
	the property \eqref{eq:semigroup_property}. 
	For any $\x_0$, we have 
\begin{align*}
\|\PPhi_{\delta} (\x_0) - \x_0\|  = \left\| \int_0^\delta \f ( \PPhi_t (\x_0) ) dt   \right\|
= \delta \eta \bigg( \frac{1}{\delta} \int_0^\delta \f ( \PPhi_t (\x_0) ) dt   \bigg),
\end{align*}
where $\eta(\x)=\| \x \|$. Since $\eta$ is a convex function, it satisfies the Jensen's inequality
$$
\eta \bigg( \frac{1}{\delta} \int_0^\delta \f ( \PPhi_t (\x_0) ) dt   \bigg) \le  \frac{1}{\delta}  \int_0^\delta \eta \big( \f ( \PPhi_t (\x_0) ) \big) dt.   
$$
Thus we obtain
\begin{align*}
\|\PPhi_{\delta} (\x_0) - \x_0\|  \le \int_0^\delta \left\|  \f ( \PPhi_t (\x_0) ) \right\| dt   
\le \delta  \sup_{t \in [0, \delta]} \| {\bf f} ( \PPhi_t (\x_0) ) \|,
\end{align*}
which implies \eqref{WKL2}. For any $\x_0 \in D_\delta$, we have $\PPhi_t (\x_0) \in D$ for $0 \le t \le \delta$. Hence 
$$
\|\PPhi_{\delta} (\x_0) - \x_0\| \le \frac{\Delta}{K}  \| {\bf f}  \|_{L^\infty(D)}, \quad \forall \x_0 \in D_\delta.
$$
This yields \eqref{WKL3}, and the proof is complete.
\end{proof}

%\begin{remark}
%		The estimate \eqref{WKL3} indicates that, when $K \gg 1$, the $\delta$-incremental flow map $\PPhi_{\delta} $ is
%	near-identity in $L^\infty$ norm. 
%	This feature implies that the flow map $\PPhi_{\delta}$ can be approximated effectively by 
%	ResNet.
%\end{remark}
This estimate can serve as a theoretical justification of the ResNet
method ($K=1$) and RT-ResNet method ($K\geq 1$). As long as $\Delta$
is reasonably small, the flow map of the underlying dynamical system
is close to identity. Therefore, it is natural to use ResNet, which
explicitly introduces the identity operator, to approximate the
``residue'' of the flow map. The norm of the DNN operator $\N$, which
approximates the residual flow map, $\PPhi_\delta-{\bf I}$, becomes
small at $O(\Delta)$. For RT-ResNet with $K>1$, its norm becomes even
smaller at $O(\Delta/K)$. We remark that it was pointed out
empirically in \cite{ChangEtAl2018} that using multiple ResNet blocks
can result in networks with smaller norm.

\subsection{Error Bound}

Let ${\bm {\mathcal N}}$ denote the neural network approximation operator to
the $\Delta$-lag flow map $\PPhi_{\Delta}$. For the proposed
ResNet \eqref{Res}, RT-ResNet \eqref{RRes}, and RS-ResNet
\eqref{Ress}, the operators can be written as
\be
\begin{split}
{\bm {\mathcal N}} & = {\bf I} + {\bf N}(\bullet ;\Theta), \quad \textrm{ResNet};\\
{\bm {\mathcal N}} & = \underbrace{\big( {\bf I} + {\bf N}(\bullet ;\Theta) \big) \circ \cdots \circ 
\big( {\bf I} + {\bf N}(\bullet ;\Theta) \big)}_{K{\rm-fold}}, \quad \textrm{RT-ResNet};\\
{\bm {\mathcal N}}  &= \big( {\bf I} + {\bf N}(\bullet ;\Theta_{K-1}) \big) \circ \cdots \circ 
\big( {\bf I} + {\bf N}(\bullet ;\Theta_{0}) \big), \quad \textrm{RS-ResNet}.
\end{split}
\ee

We now derive a general error bound for the solution approximation
using the DNN operator ${\bm {\mathcal N}}$.  This bound serves a
general guideline for the error growth. More specific error
bounds for each different network structure are more involved and will
be pursued in a future work.

Let ${\y}^{(m)}$ denote the solution of the approximate model at time
$t^{(m)}:=t_0+m\Delta$. Let ${\mathcal E}^{(m)}:= \| {\y}^{(m)} - \x (
t^{(m)} ) \|$ denote the error,  
$j=0,1,\dots,m$.

\begin{theorem}\label{thm:err}
Assume that the same assumptions in Lemma \ref{lem1} hold, and let us
further assume
\begin{enumerate}
	\item $\big\| {\bm {\mathcal N}}- \PPhi_{\Delta} \big\|_{L^\infty (D_\Delta)} < +\infty$,
	\item ${\y}^{(i)}, \x(t^{(i)}) \in D_\Delta$ for $0\le i \le m-1$,
\end{enumerate}
then we have
\begin{equation}\label{eq:WKLest2}
{\mathcal E}^{(m)} \le \big( 1 + {\rm e}^{L \Delta}   \big)^m {\mathcal E}^{(0)} 
+ \big\| {\bm {\mathcal N}}- \PPhi_{\Delta} \big\|_{L^\infty (D_\Delta)}  \frac{  \big( 1 + {\rm e}^{L \Delta}   \big)^m - 1 }{ {\rm e}^{L \Delta}  }.
\end{equation}
\end{theorem}

\begin{proof}
The triangle inequality implies that 
\begin{align*}
{\mathcal E}^{(m)} & = \big\| \y^{(m-1)} + {\bm {\mathcal N}} (  \y^{(m-1)} ) 
- \x ( t^{(m-1)} ) - \PPhi_{\Delta} ( \x ( t^{(m-1)} )  )  \big \|
\\
&\le \big\| \y^{(m-1)}  - \x ( t^{(m-1)} )   \big \| + \big \| {\bm {\mathcal N}} (  \y^{(m-1)} )  - \PPhi_{\Delta} ( \x ( t^{(m-1)} )  )  \big\| 
\\
& \le \big\| \y^{(m-1)}  - \x ( t^{(m-1)} )   \big \| + \big \| {\bm {\mathcal N}} (  \y^{(m-1)} )  - 
\PPhi_{\Delta} ( \y ^{(m-1)}  )  \big\| 
\\
& \qquad 
+ \big \|  \PPhi_{\Delta} ( \y ^{(m-1)}  )  - \PPhi_{\Delta} ( \x ( t^{(m-1)} )  )  \big\| 
\\
& \le  \big\| \y^{(m-1)}  - \x ( t^{(m-1)} )   \big \| + \big \| {\bm {\mathcal N}} - 
\PPhi_{\Delta}  \big\|_{L^\infty(D_\Delta)} 
\\
& \qquad 
+ {\rm e}^{L \Delta} \big \|   \y ^{(m-1)}    -  \x ( t^{(m-1)} )   \big\|
\\
& = \big( 1+ {\rm e}^{L \Delta} \big)  {\mathcal E}^{(m-1)}  + \big \| {\bm {\mathcal N}} - 
\PPhi_{\Delta}  \big\|_{L^\infty(D_\Delta)} ,
\end{align*}
where the Lipschitz continuity of the flow map, shown in \eqref{eq:WKLest1}, has been used in the last inequality. Recursively using the above estimate gives
\begin{align*}
{\mathcal E}^{(m)} & 
\le \big( 1+ {\rm e}^{L \Delta} \big)   {\mathcal E}^{(m-1)}  + \big \| {\bm {\mathcal N}} - 
\PPhi_{\Delta}  \big\|_{L^\infty(D_\Delta)} 
\\
&\le \big( 1+ {\rm e}^{L \Delta} \big)^2   
  {\mathcal E}^{(m-2)}  
  + \big \| {\bm {\mathcal N}} - 
\PPhi_{\Delta}  \big\|_{L^\infty(D_\Delta)} \Big( 1 + \big( 1+ {\rm e}^{L \Delta} \big) \Big)
\\
&
\le \cdots
\\
&\le  \big( 1+ {\rm e}^{L \Delta} \big)^m   
{\mathcal E}^{(0)}  
+ \big \| {\bm {\mathcal N}} - 
\PPhi_{\Delta}  \big\|_{L^\infty(D_\Delta)} \sum_{i=0}^{m-1} \big( 1+ {\rm e}^{L \Delta} \big) ^i.
\end{align*}	
The proof is complete.	
\end{proof}

%Analysis of RS-ResNet is more involved, as its
%adaptive nature makes error bound dependent on the data. This will be
%pursued in a separate work.

\section{Numerical Examples} \label{sec:examples}

In this section we present numerical examples to verify the properties
of the proposed methods. In all the examples, we generate the training
data pairs
$\{(\z_j^{(1)}, \z_j^{(2)}\}_{j=1}^J$ in the following way:
\begin{itemize}
	\item Generate $J$ points $\{\z_j^{(1)}\}_{j=1}^J$ from
          uniform distribution over a computational domain $D$. The
          domain $D$ is a region in which we are interested in the
          solution behavior. It is typically chosen to be a hypercube
          prior to the computation.
	\item For each $j$, starting from $\z_j^{(1)}$, we march
          forward for a time lag $\Delta$ the underlying governing
          equation, using a highly accurate
		standard ODE solver, to generate
                $\z_j^{(2)}$. In our examples we set
                $\Delta=0.1$. 
\end{itemize}
We remark that the time lag $\Delta=0.1$ is relatively coarse and
prevents accurate estimate of time derivatives via numerical
differentiation. Since our proposed methods employ the integral form
of the underlying equation, this difficulty is circumvented. 
The random sampling of the solution trajectories of length $\Delta$
follows from the work of \cite{WuXiu_JCPEQ18}, where it was
established that such kind of dense sampling of short trajectories is
highly effective for equation recovery.

All of our network models, ResNet, RT-ResNet, and RS-ResNet, are
trained via the loss function \eqref{loss} and by using the
open-source Tensorflow library \cite{tensorflow2015}.
%We train the model by minimizing the mean squared loss function \eqref{loss}. The optimization problem is solved with the Adam algorithm \cite{Adam}.
The training data set is divided into mini-batches of size $10$. And
we typically train the model for
$500$ epochs and reshuffle the training data in each epoch.
All the weights are initialized randomly from Gaussian
distributions and all the biases are initialized to be zeros. 
%All the model
%construction, training and prediction are performed with the
%open-source Tensorflow library \cite{tensorflow2015}. 

After training the network models satisfactorily, using the data of
$\Delta=0.1$ time lag, we march the trained
network models further forward in time and compare the results against the
reference states, which are produced by high-order numerical solvers
of the true
underlying governing equations. We march the trained network systems
up to $t\gg \Delta$ to examine their (relatively) long-term
behaviors. For the two linear examples, we set $t=2$; and for the two
nonlinear examples, we set $t=20$.

\subsection{Linear ODEs}

We first study two linear ODE systems, as textbook examples. In both
examples, our one-step ResNet method has 3 hidden layers, each of
which has 30 neurons. For the multi-step RT-ResNet and RS-ResNet
methods, they both have $3$ ResNet blocks ($K=3$), each of which
contains 3 hidden layers with 20 neurons in each layer.

\subsubsection*{Example 1}
We first consider the following two-dimensional linear ODE with $\mathbf{x}=(x_1, x_2)$
\begin{equation}
	\label{eq:example1}
	\begin{cases}
		\dot{x}_1=x_1+x_2-2,\\
		\dot{x}_2=x_1-x_2.
	\end{cases}
\end{equation}
The computational domain $D$ is taken to be $D=[0, 2]^2$. 

\begin{figure}[!htb]
	\centering
	\begin{subfigure}[b]{0.32\textwidth}
	\begin{center}
		\includegraphics[width=1.0\linewidth]{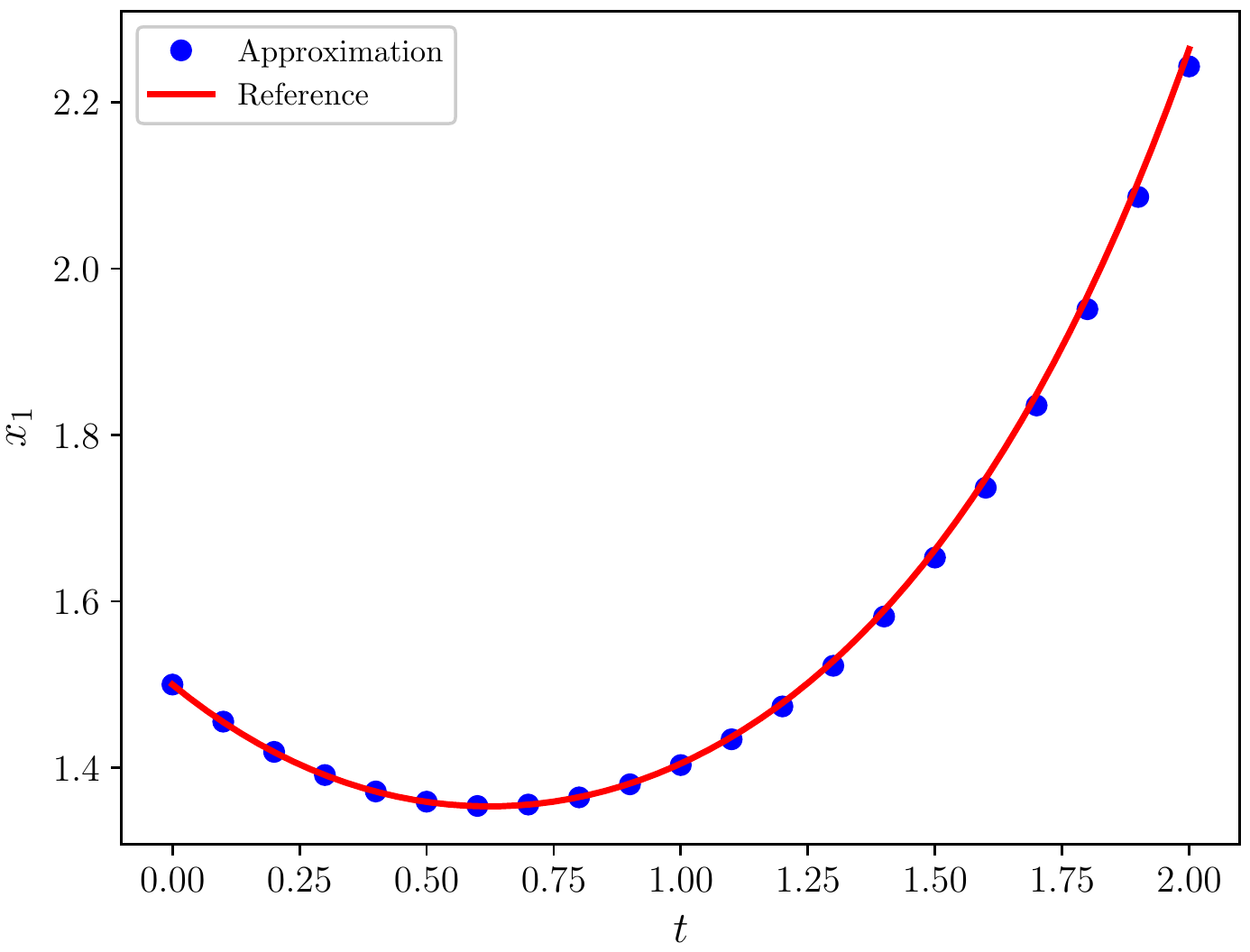}
		\caption{$x_1$, ResNet}
	\end{center}
	\end{subfigure}
	\begin{subfigure}[b]{0.32\textwidth}
		\centering
		\includegraphics[width=1.0\linewidth]{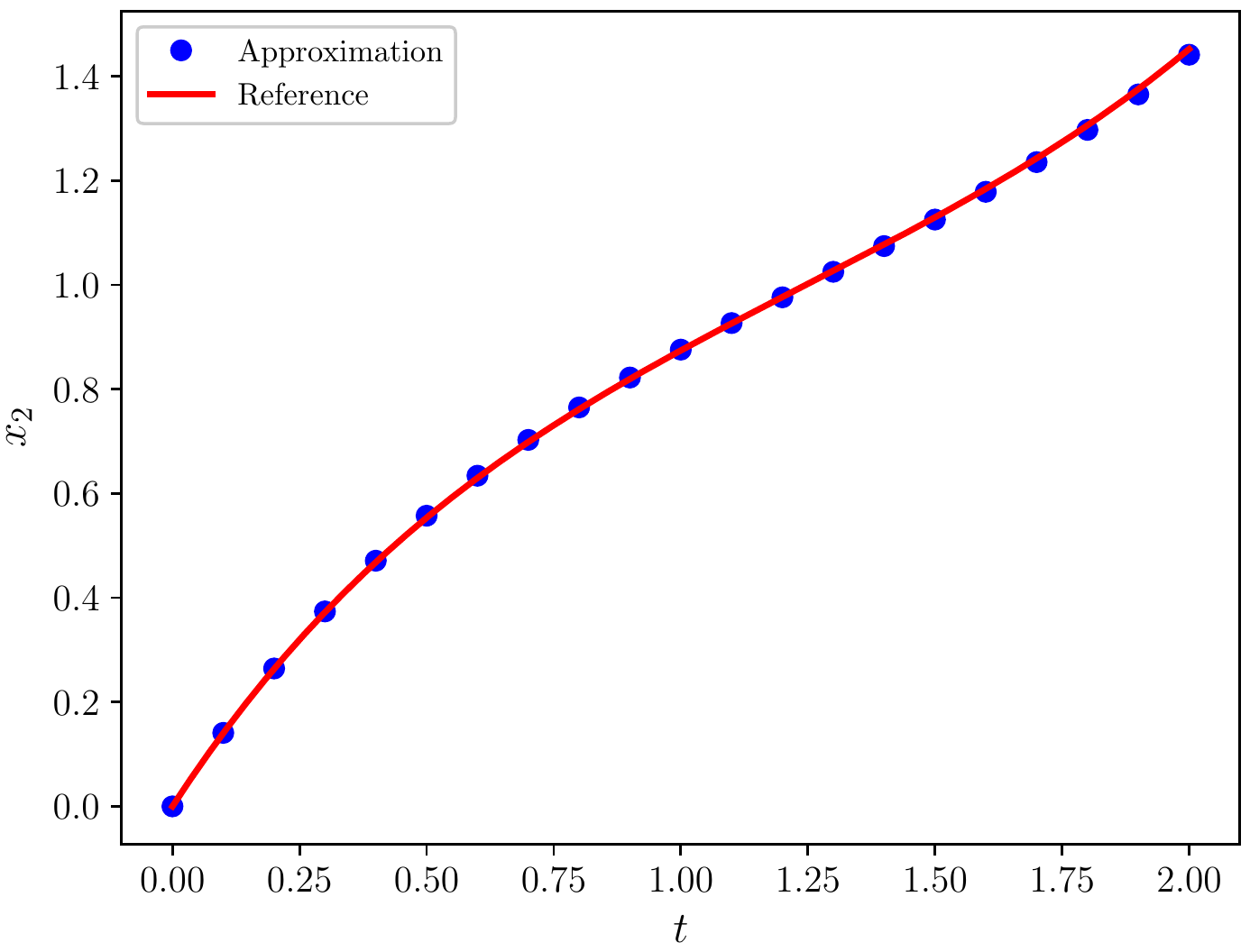}
		\caption{$x_2$, ResNet}
	\end{subfigure}
	\begin{subfigure}[b]{0.32\textwidth}
	\begin{center}
		\includegraphics[width=1.0\linewidth]{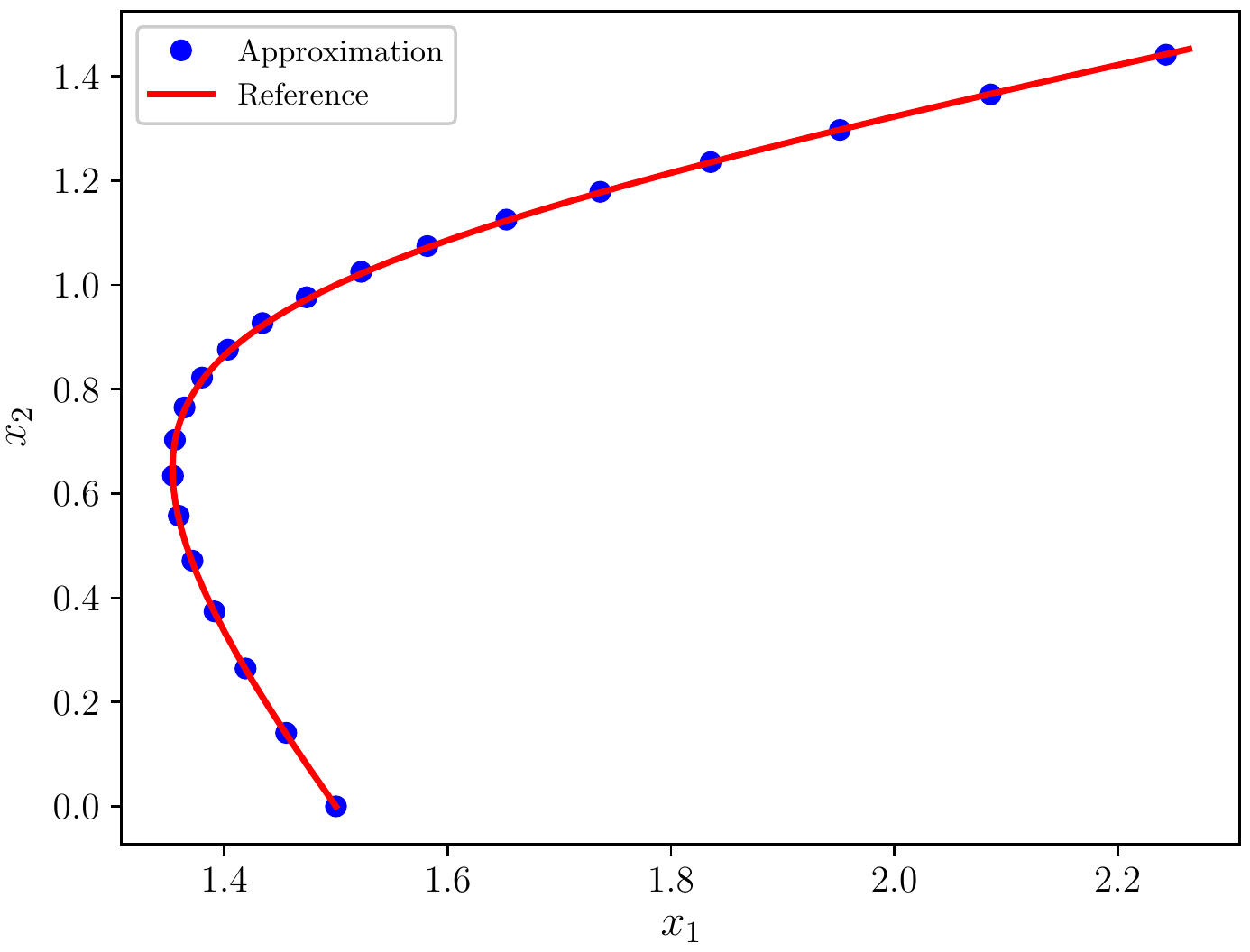}
		\caption{phase plot, ResNet}
	\end{center}
	\end{subfigure}
	\begin{subfigure}[b]{0.32\textwidth}
	\begin{center}
		\includegraphics[width=1.0\linewidth]{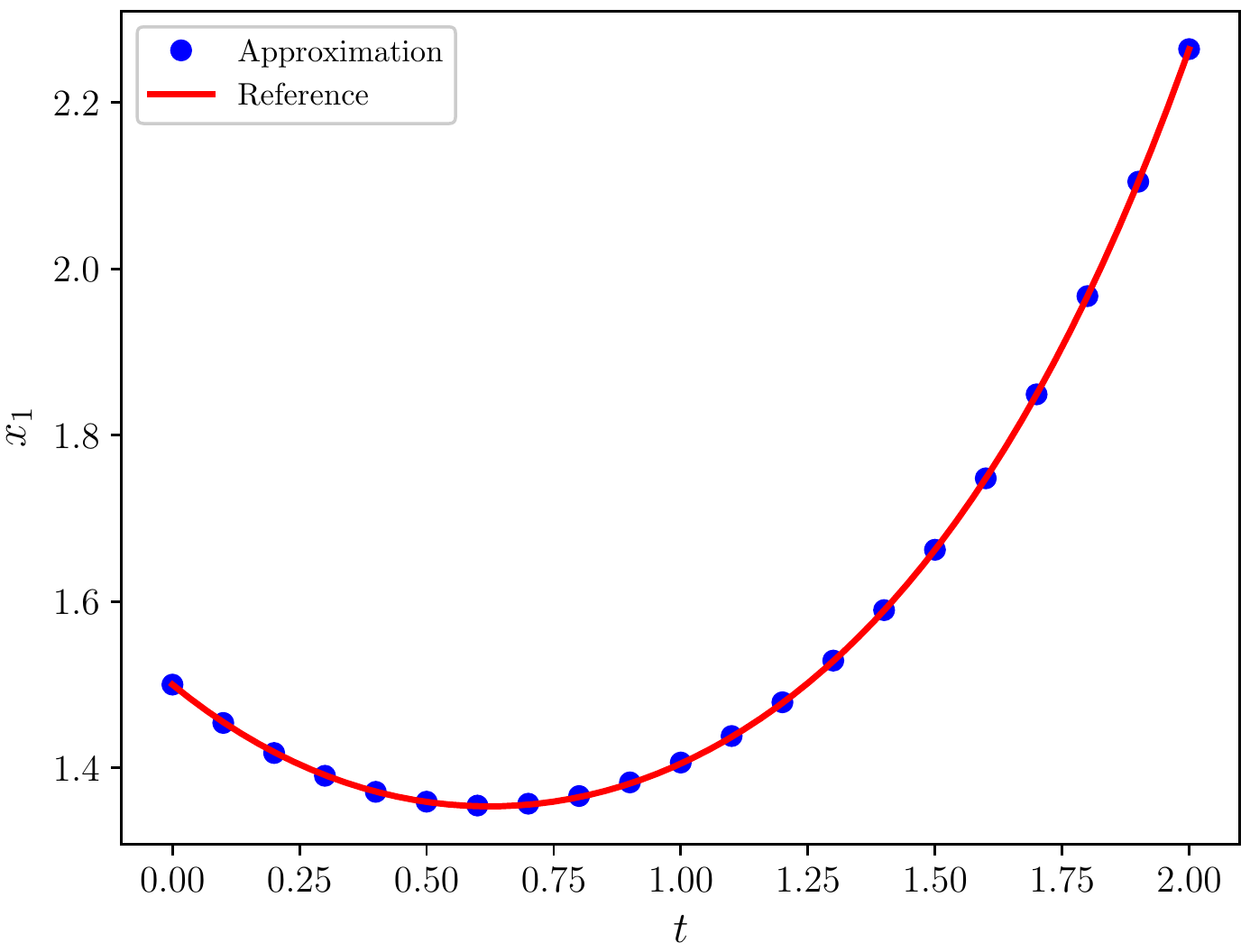}
		\caption{$x_1$, RT-ResNet}
	\end{center}
	\end{subfigure}
	\begin{subfigure}[b]{0.32\textwidth}
		\centering
		\includegraphics[width=1.0\linewidth]{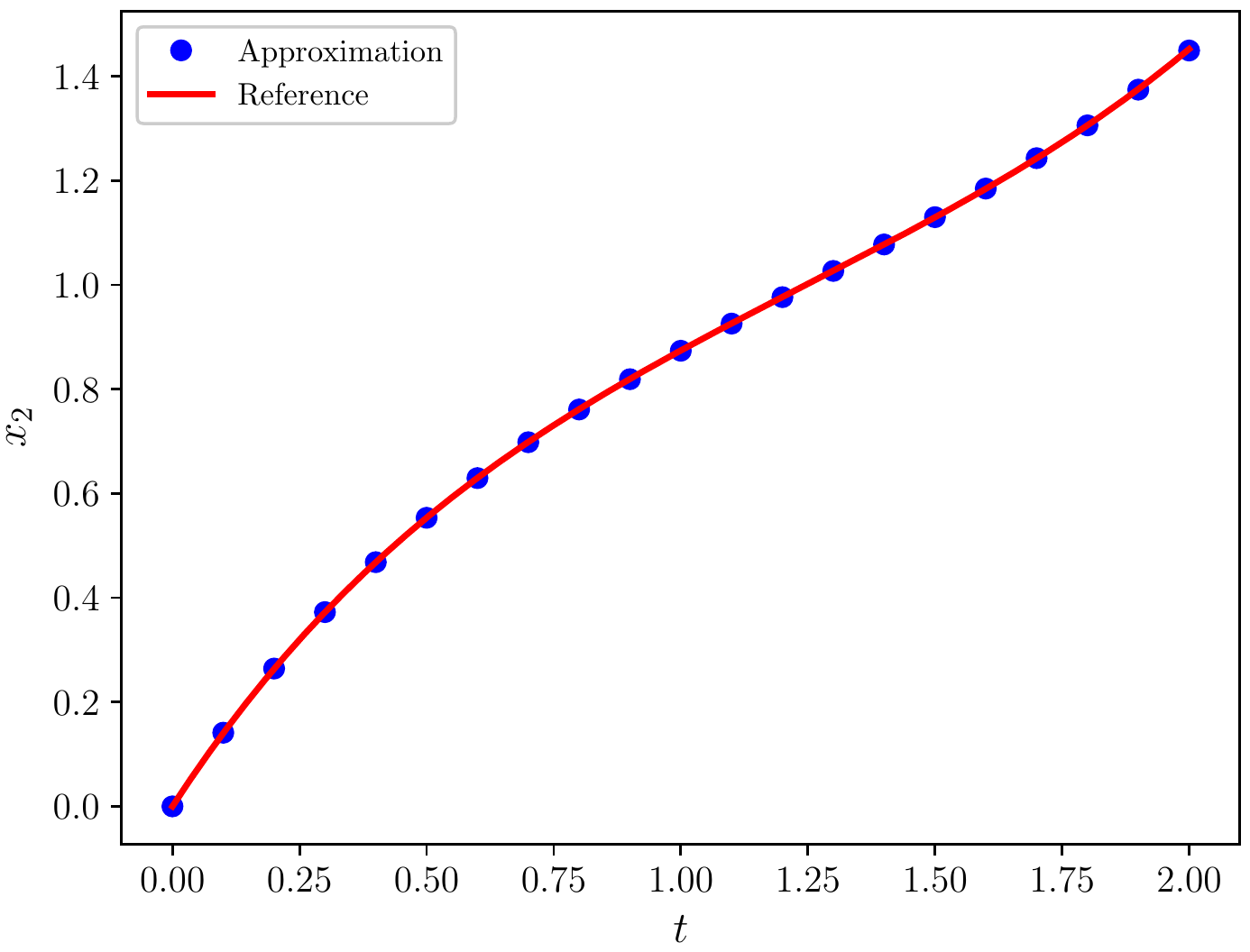}
		\caption{$x_2$, RT-ResNet}
	\end{subfigure}
	\begin{subfigure}[b]{0.32\textwidth}
	\begin{center}
		\includegraphics[width=1.0\linewidth]{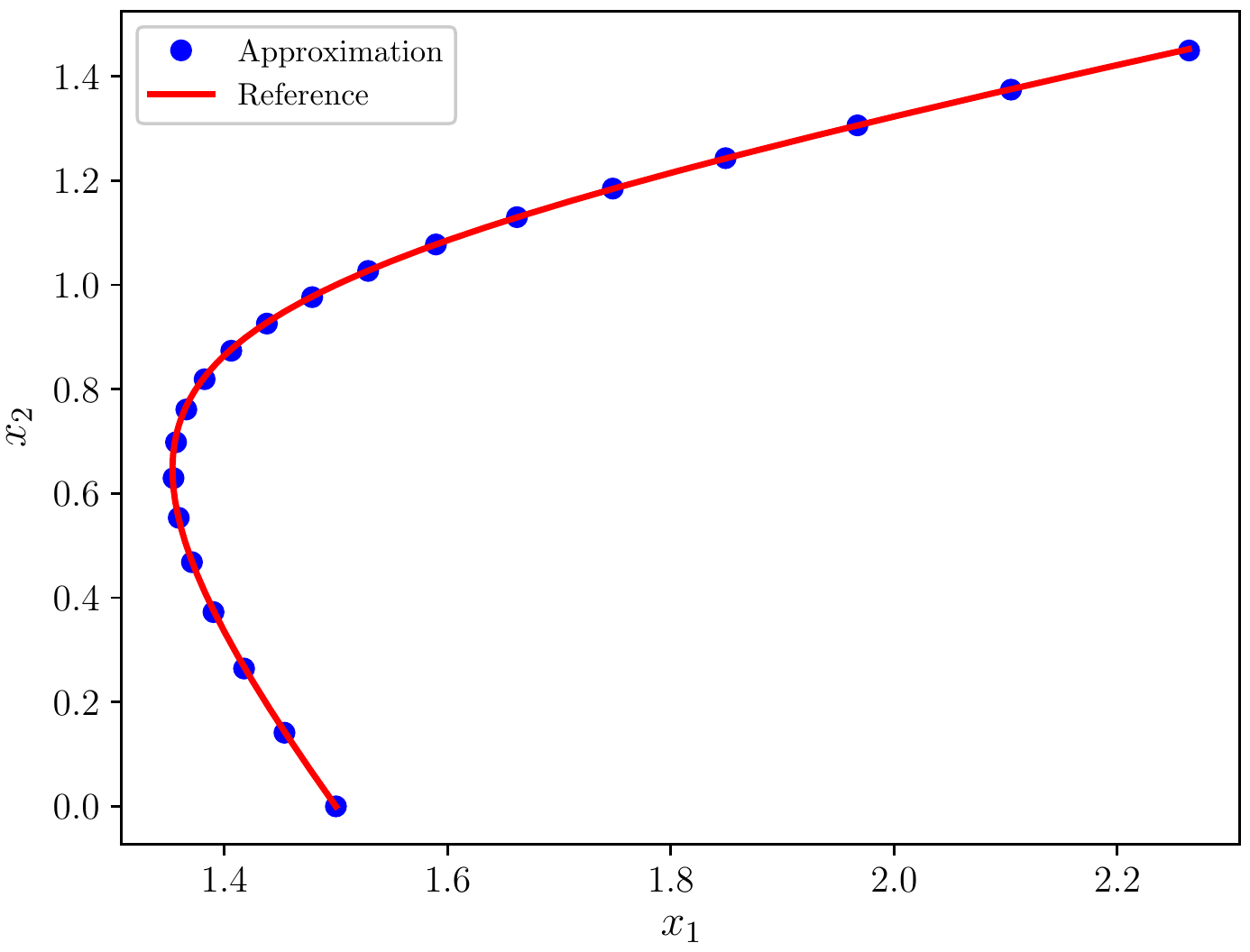}
		\caption{phase plot, RT-ResNet}
	\end{center}
	\end{subfigure}
	\begin{subfigure}[b]{0.32\textwidth}
	\begin{center}
		\includegraphics[width=1.0\linewidth]{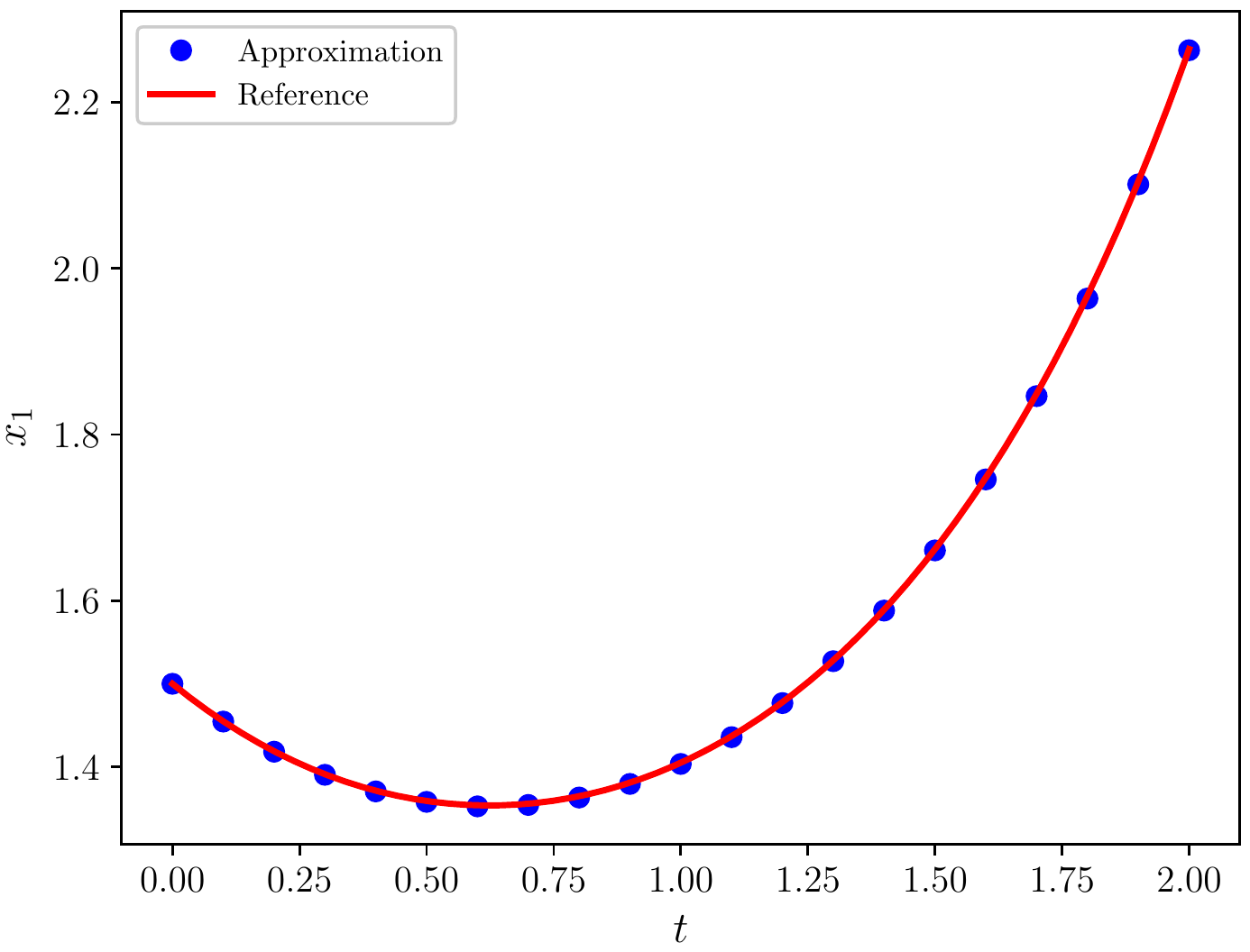}
		\caption{$x_1$, RS-ResNet}
	\end{center}
	\end{subfigure}
	\begin{subfigure}[b]{0.32\textwidth}
		\centering
		\includegraphics[width=1.0\linewidth]{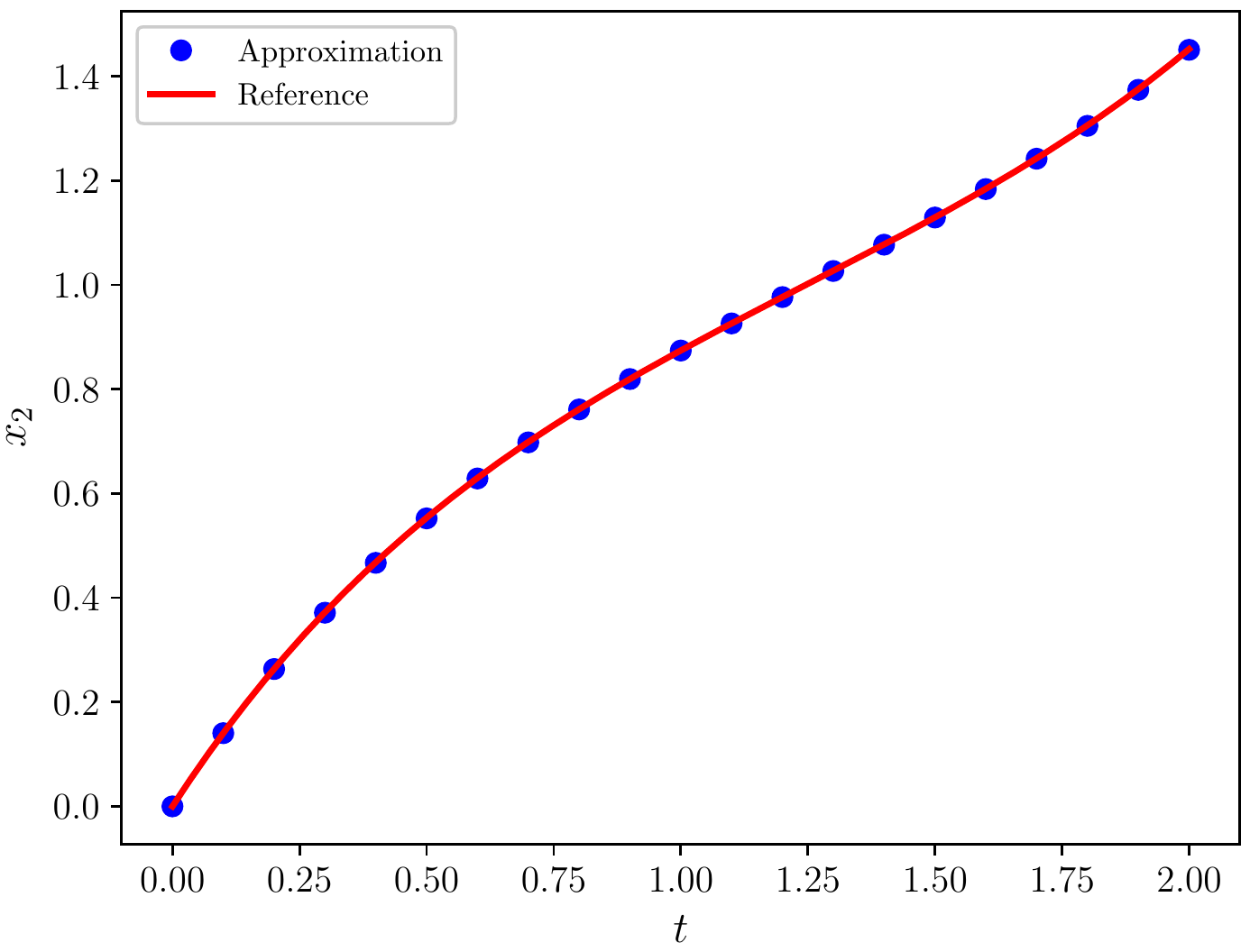}
		\caption{$x_2$, RS-ResNet}
	\end{subfigure}
	\begin{subfigure}[b]{0.32\textwidth}
	\begin{center}
		\includegraphics[width=1.0\linewidth]{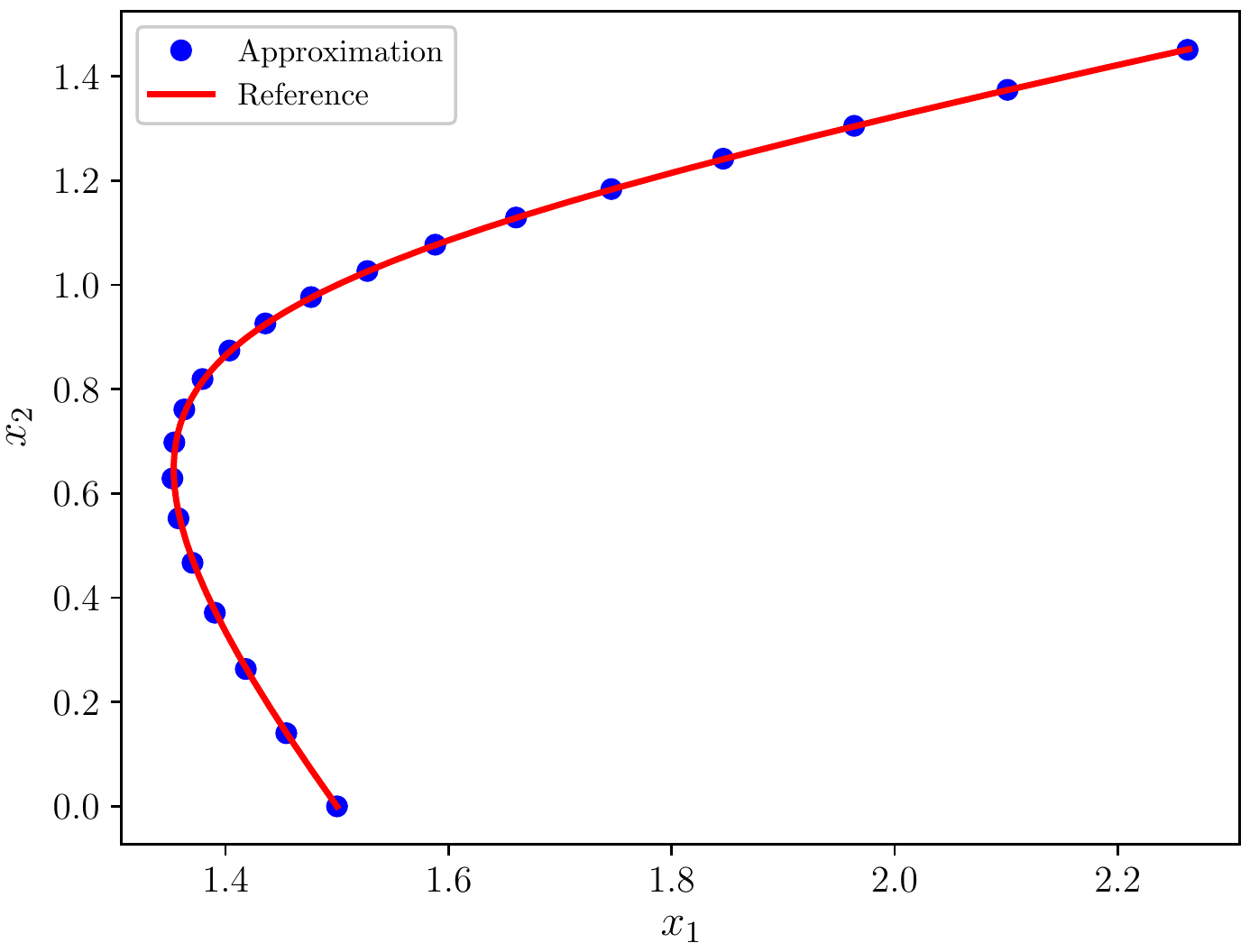}
		\caption{phase plot, RS-ResNet}
	\end{center}
	\end{subfigure}
	\caption{Trajectory and phase plots for the Example 1 with
          $\x_0=(1.5, 0)$ for $t\in [0,2]$.
Top row: one-step ResNet model;
Middle row: Multi-step RT-ResNet model;
Bottom row: Multi-step RS-ResNet model.
%	From top to bottom: approximations by the single-step ResNet,
%	the RT-ResNet and the RS-ResNet. Circles represent the
%	approximations at discrete time points. Lines represent
%	reference solutions.
}
	\label{fig:ex1}
\end{figure}

Upon training the three network models satisfactorily, using the
$\Delta=0.1$ data pairs, we march the
trained models further in time up to $t=2$. The results are shown in
Figure \ref{fig:ex1}. We observe that all three network models produce
accurate prediction results for time up to $t=2$.
%We test all the three DNN models: ResNet, RT-ResNet and RS-ResNet, which are proposed in Section \ref{sec:method}. The ResNet contains $3$ hidden layers and each layer contains $30$ nodes. As for the other two models, they both have three ResNet block and
%each ResNet block contains $3$ hidden layers with 
%$20$ nodes in each layer. The trained DNN approximates the flow map with $\Delta=0.1$ lag. We then compose it with itself for $20$ times to generate
%the approximate trajectories and the phase plots from $t=0$ to $t=2$. In Figure \ref{fig:ex1} we present the results generated by these three models. 
%We see that all of these three models deliver very good approximation. 

\begin{figure}[!htb]
	\centering
	\begin{subfigure}[b]{0.32\textwidth}
	\begin{center}
		\includegraphics[width=0.9\linewidth]{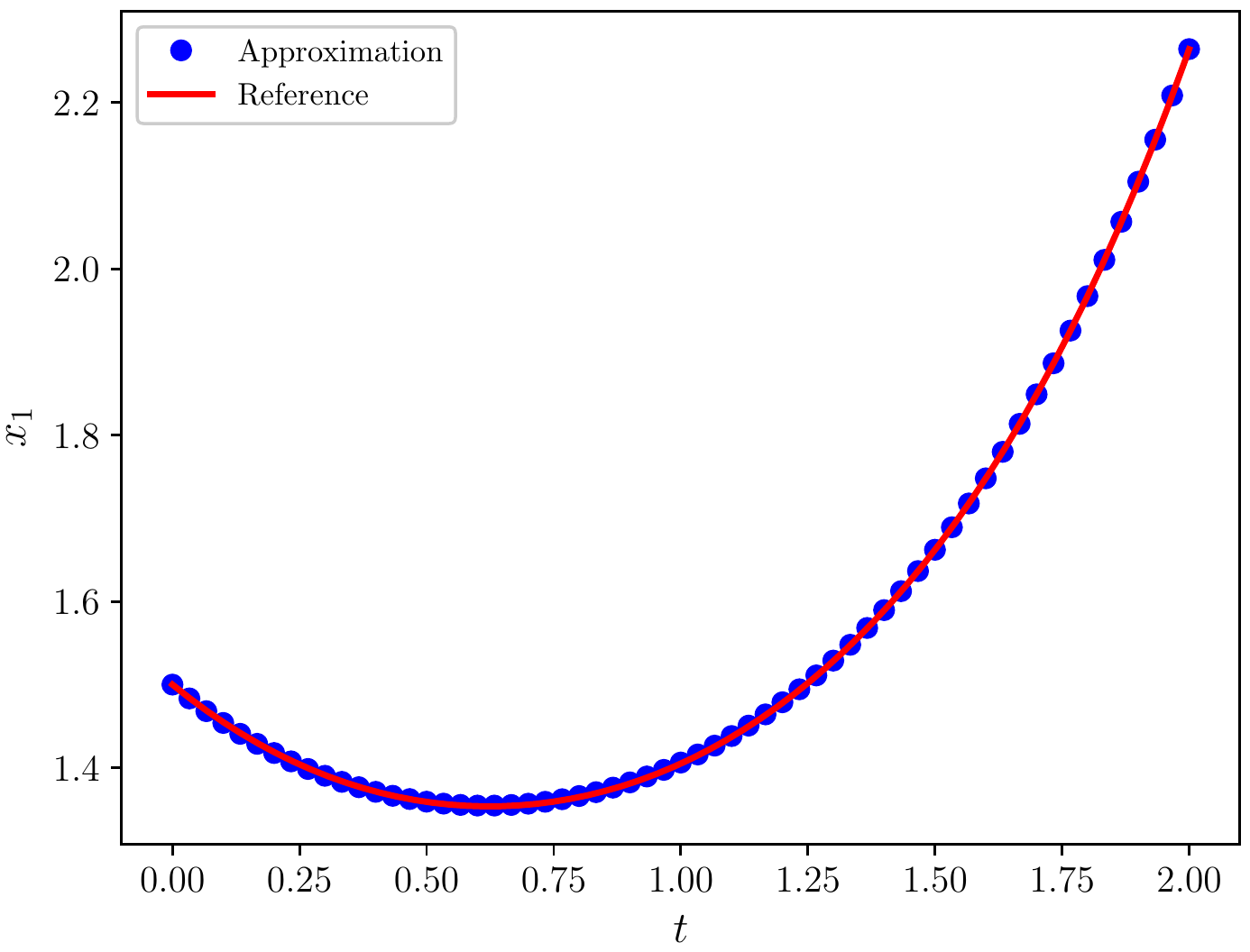}
		\caption{$x_1$}
	\end{center}
	\end{subfigure}
	\begin{subfigure}[b]{0.32\textwidth}
		\centering
		\includegraphics[width=1.0\linewidth]{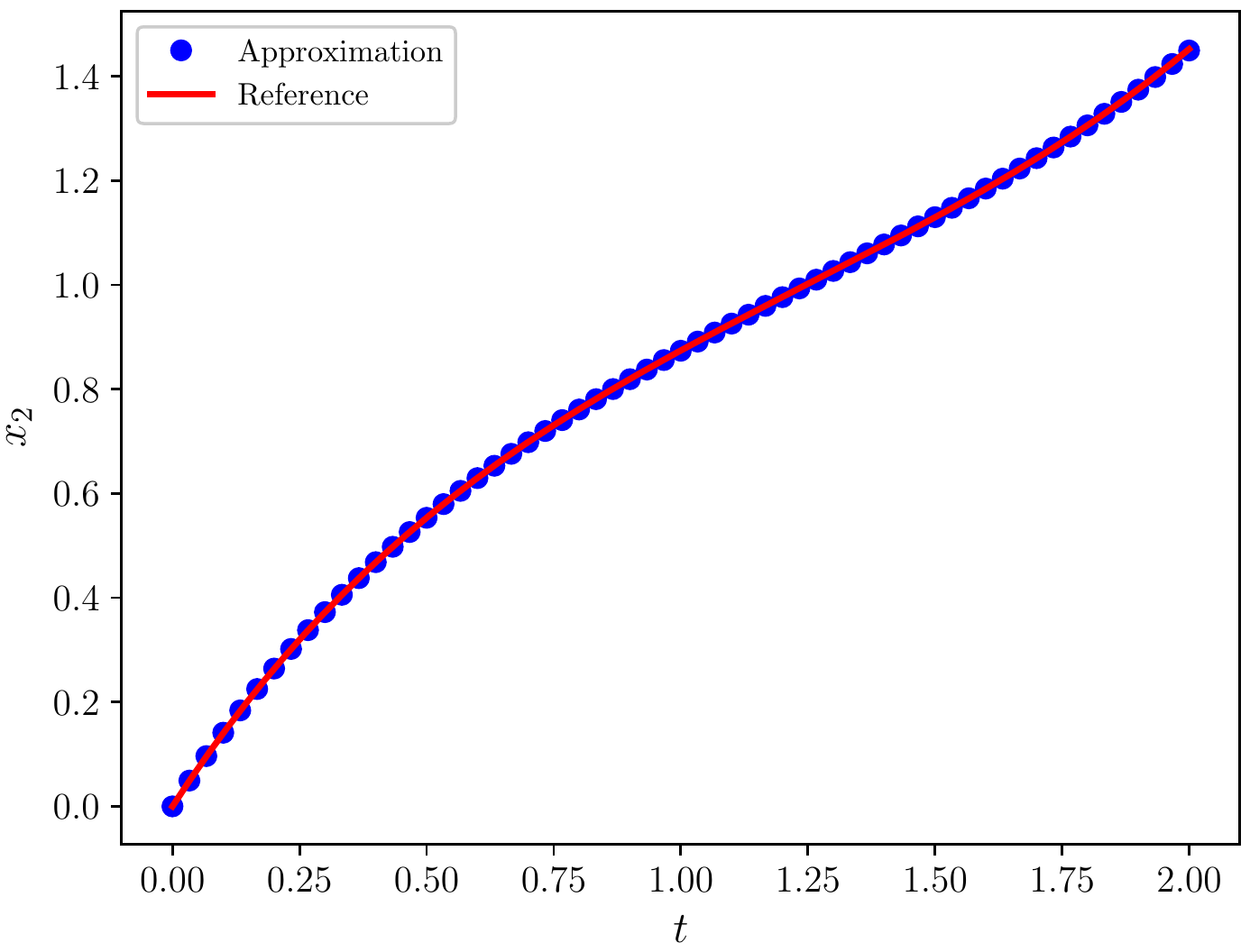}
		\caption{$x_2$}
	\end{subfigure}
	\begin{subfigure}[b]{0.32\textwidth}
	\begin{center}
		\includegraphics[width=1.0\linewidth]{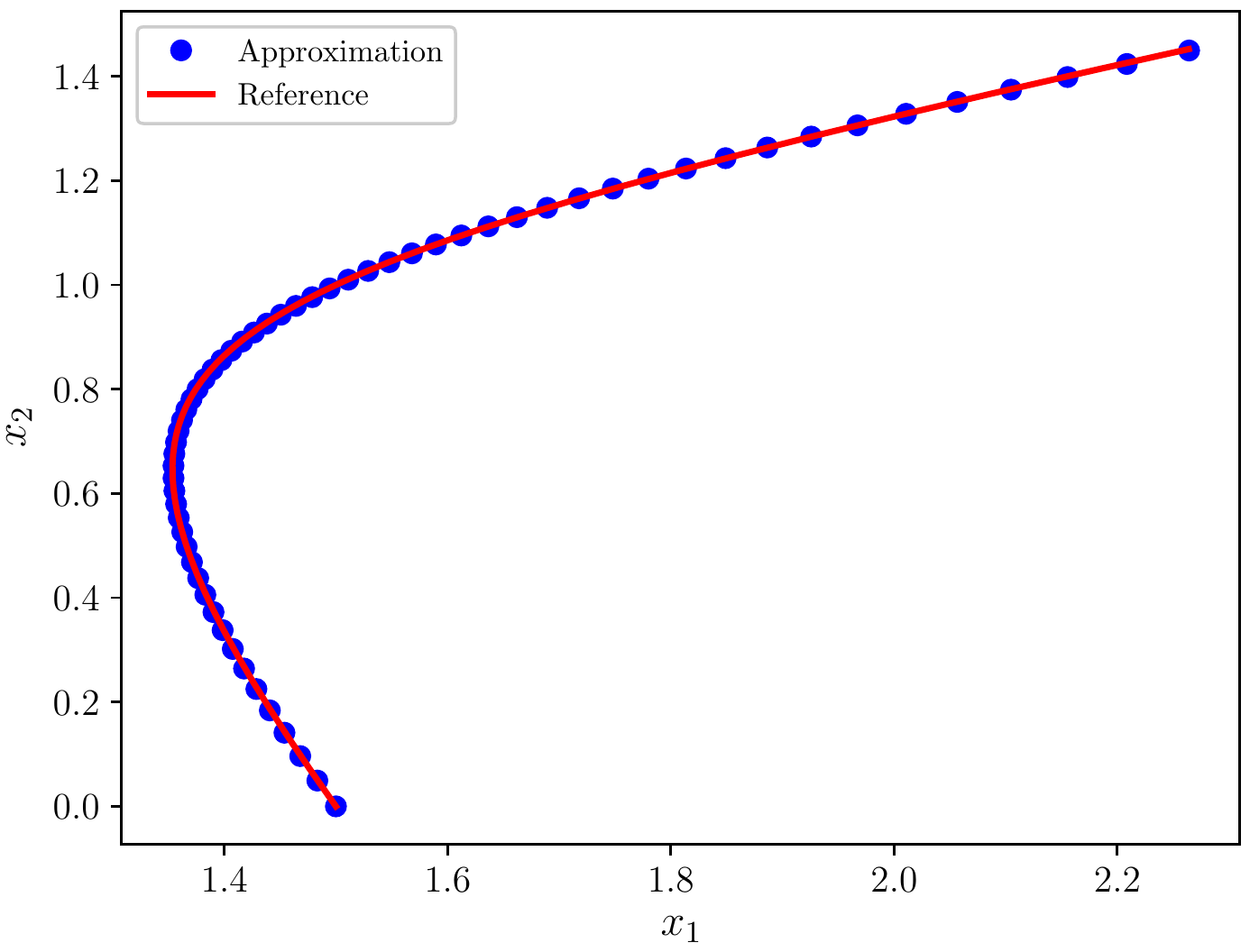}
		\caption{phase plot}
	\end{center}
	\end{subfigure}
	\caption{Trajectory and phase plots for Example 1 with
          $\x_0=(1.5, 0)$ using RT-ResNet model with $K=3$. The
          solutions are marched into over time step $\delta = \Delta/K
        = 1/30$.}
	\label{fig:ex1_delta}
\end{figure}

As discussed in Section \ref{sec:multi-step_recurrent}, the multi-step RT-ResNet method
is able to produce an approximation over a smaller
time step $\delta=\Delta/K$, which in this case is $\delta = 1/30$
(with $K=3$). The trained RT-ResNet model then allows us to produce predictions over
the finer time step $\delta$.
%Moreover, for the RT-ResNet model, as discussed in Section \ref{sec:multi-step_recurrent},
%each single ResNet block $\N(\x;\Theta)$ is approximating the $\delta$-effective
%increment $\pphi_\delta(\x;\f)$. This means that, we can potentially generate more detailed
%approximations by composing the single ResNet block with itself, if the DNN is well trained. 
%
In Figure \ref{fig:ex1_delta}, we show the time marching of the
trained RT-ResNet model for up to $t=2$ using the smaller time step
$\delta$. The results again agree very well with the reference
solution. This demonstrates the capability of RT-ResNet -- it allows
us to produce accurate predictions with a resolution higher than that
of the 
given data, i.e., $\delta < \Delta$. On the other hand, our numerical tests
also reveal that the training of RT-ResNet with $K>1$ becomes more
involving -- more training data are typically required and convergence
can be slower, compared to the training of the one-step ResNet
method. Similar behavior is also observed in multi-step RS-ResNet
method with $K>1$.
%RT-ResNet, the model is able to recover the dynamics on the finer time grid with $\delta=1/30$. But we want to
%emphasize, even though, given the data pairs over $\Delta$ time
%interval, the ResNet is potentially possible to recover the
%$\delta$-increment, the training usually takes much more epochs than
%merely learning the $\Delta$-increment. 
The development of 
efficient training procedures for multi-step RT-ResNet and RS-ResNet
methods is necessary and will be pursued in a future work.
%%%%%%%%%

\subsubsection*{Example 2}

\begin{figure}[!htb]
	\centering
	\begin{subfigure}[b]{0.32\textwidth}
	\begin{center}
		\includegraphics[width=1.0\linewidth]{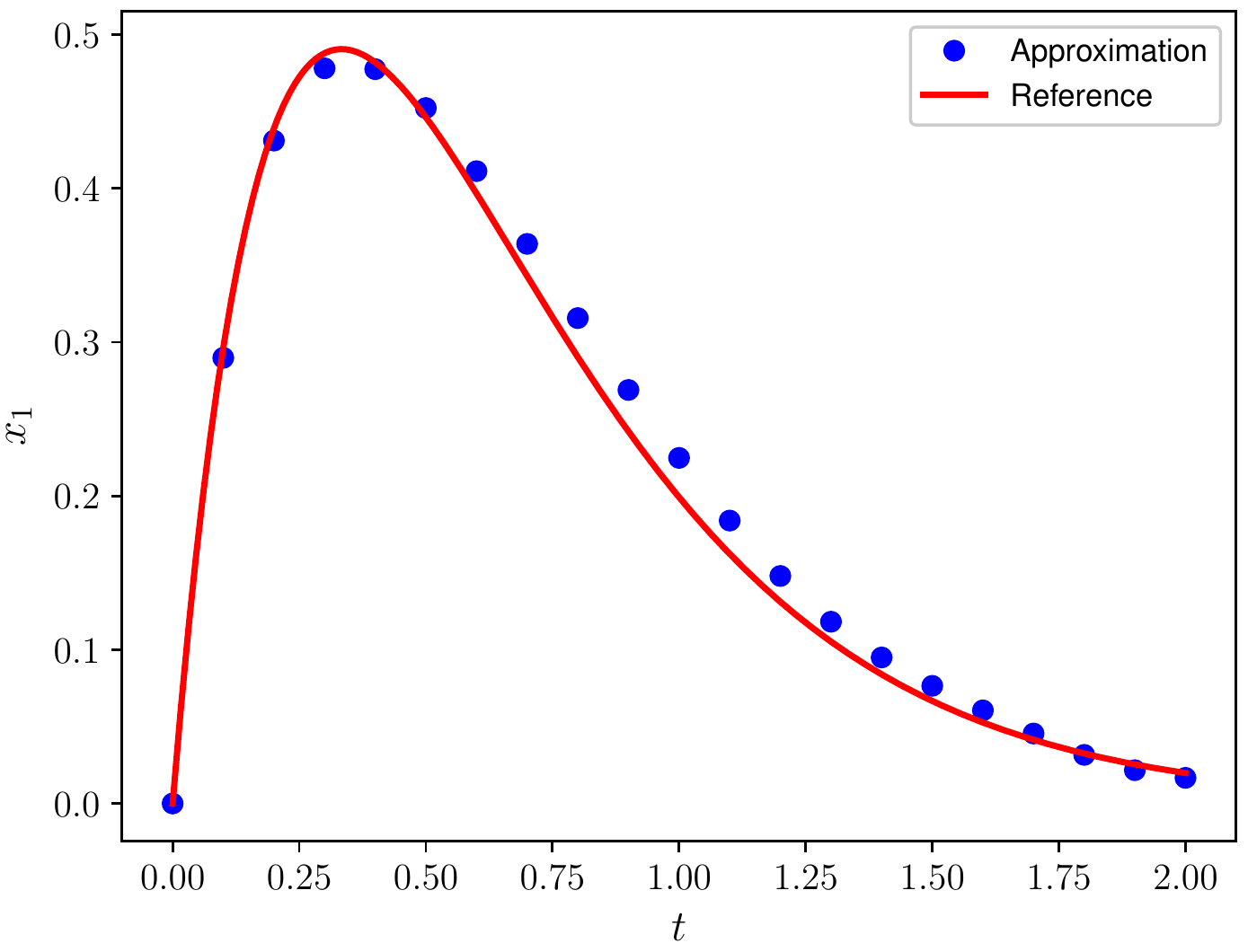}
		\caption{$x_1$, ResNet}
	\end{center}
	\end{subfigure}
	\begin{subfigure}[b]{0.32\textwidth}
		\centering
		\includegraphics[width=1.0\linewidth]{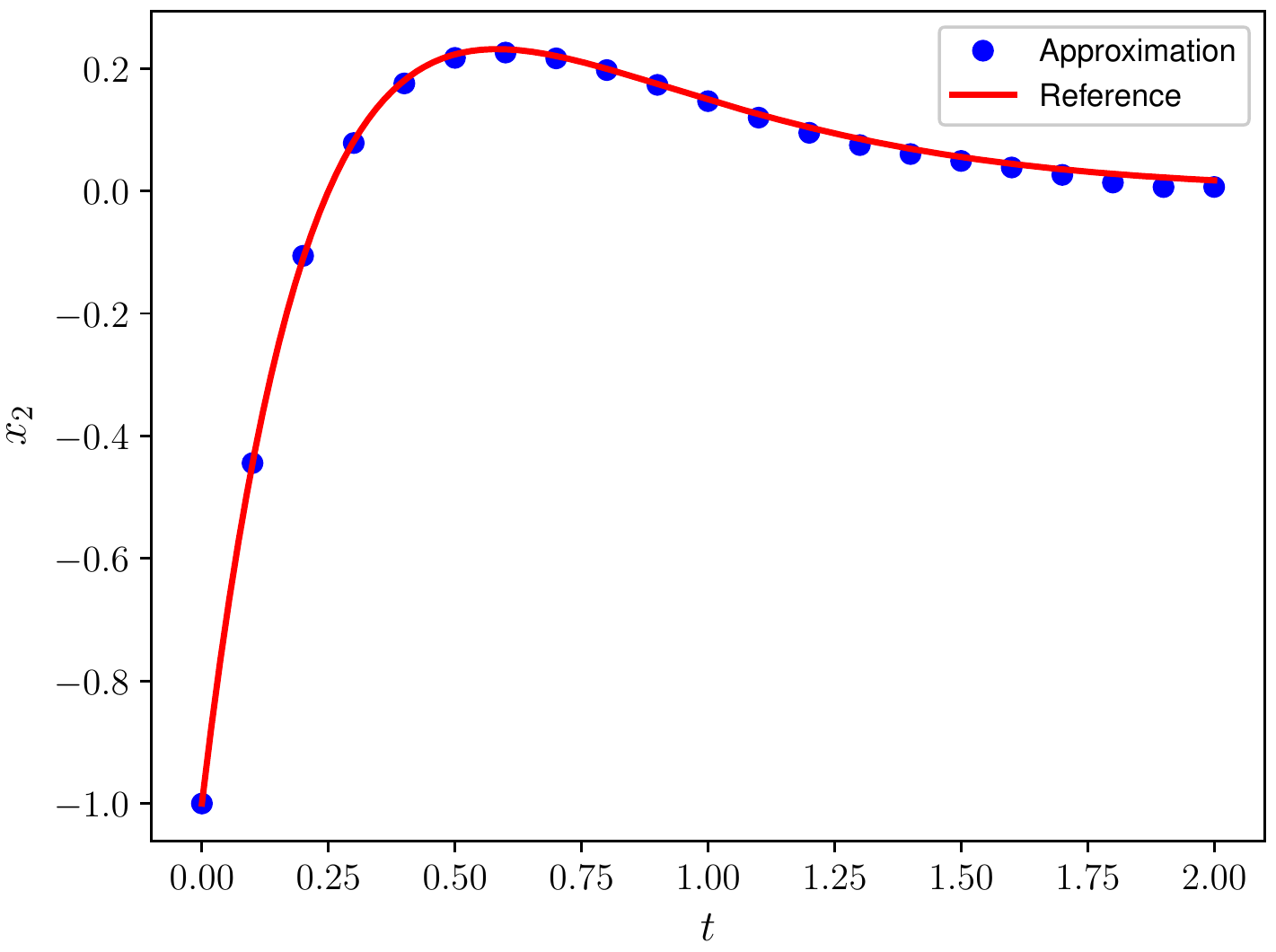}
		\caption{$x_2$, ResNet}
	\end{subfigure}
	\begin{subfigure}[b]{0.32\textwidth}
	\begin{center}
		\includegraphics[width=1.0\linewidth]{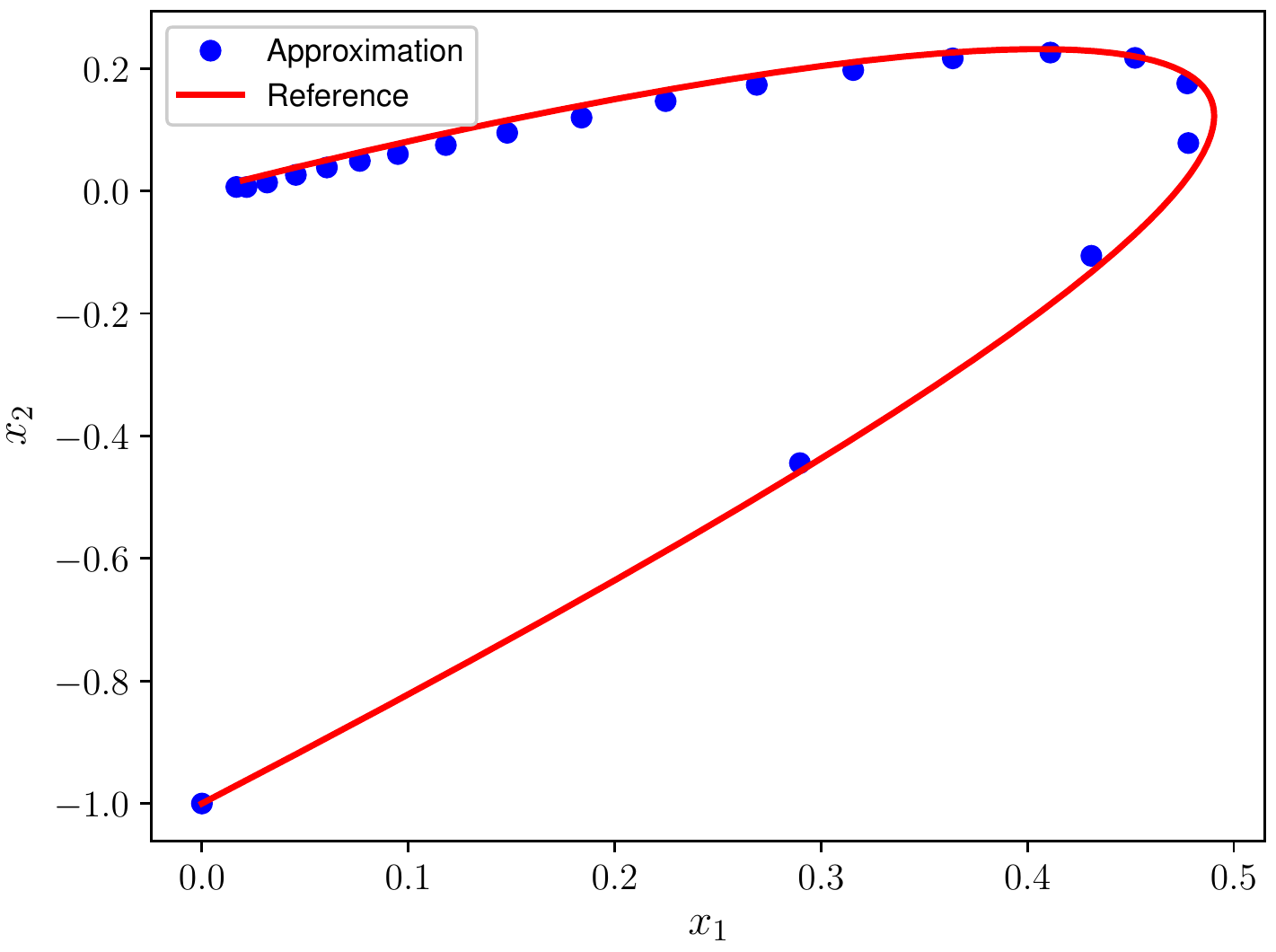}
		\caption{phase plot, ResNet}
	\end{center}
	\end{subfigure}
	\begin{subfigure}[b]{0.32\textwidth}
	\begin{center}
		\includegraphics[width=1.0\linewidth]{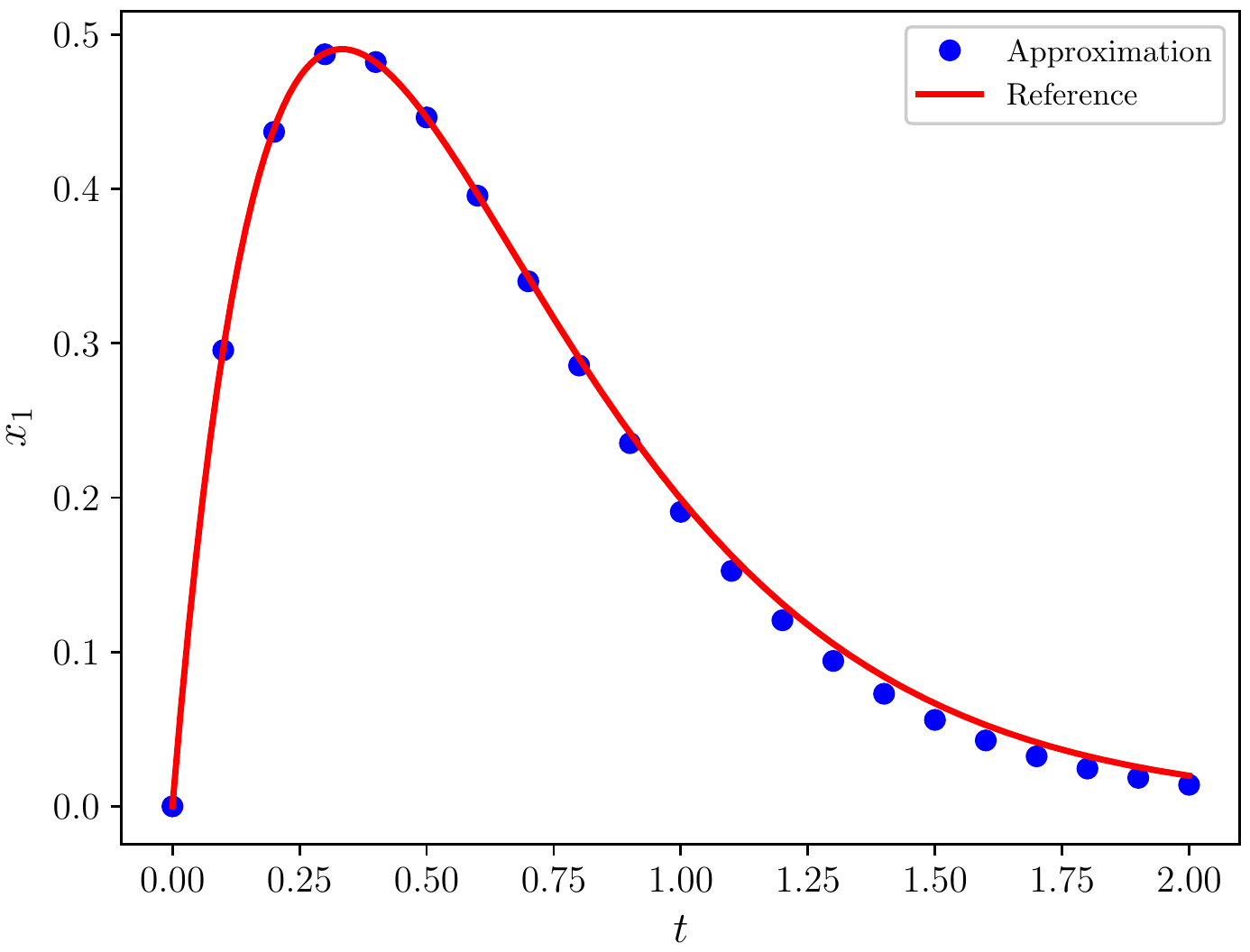}
		\caption{$x_1$, RT-ResNet}
	\end{center}
	\end{subfigure}
	\begin{subfigure}[b]{0.32\textwidth}
		\centering
		\includegraphics[width=1.0\linewidth]{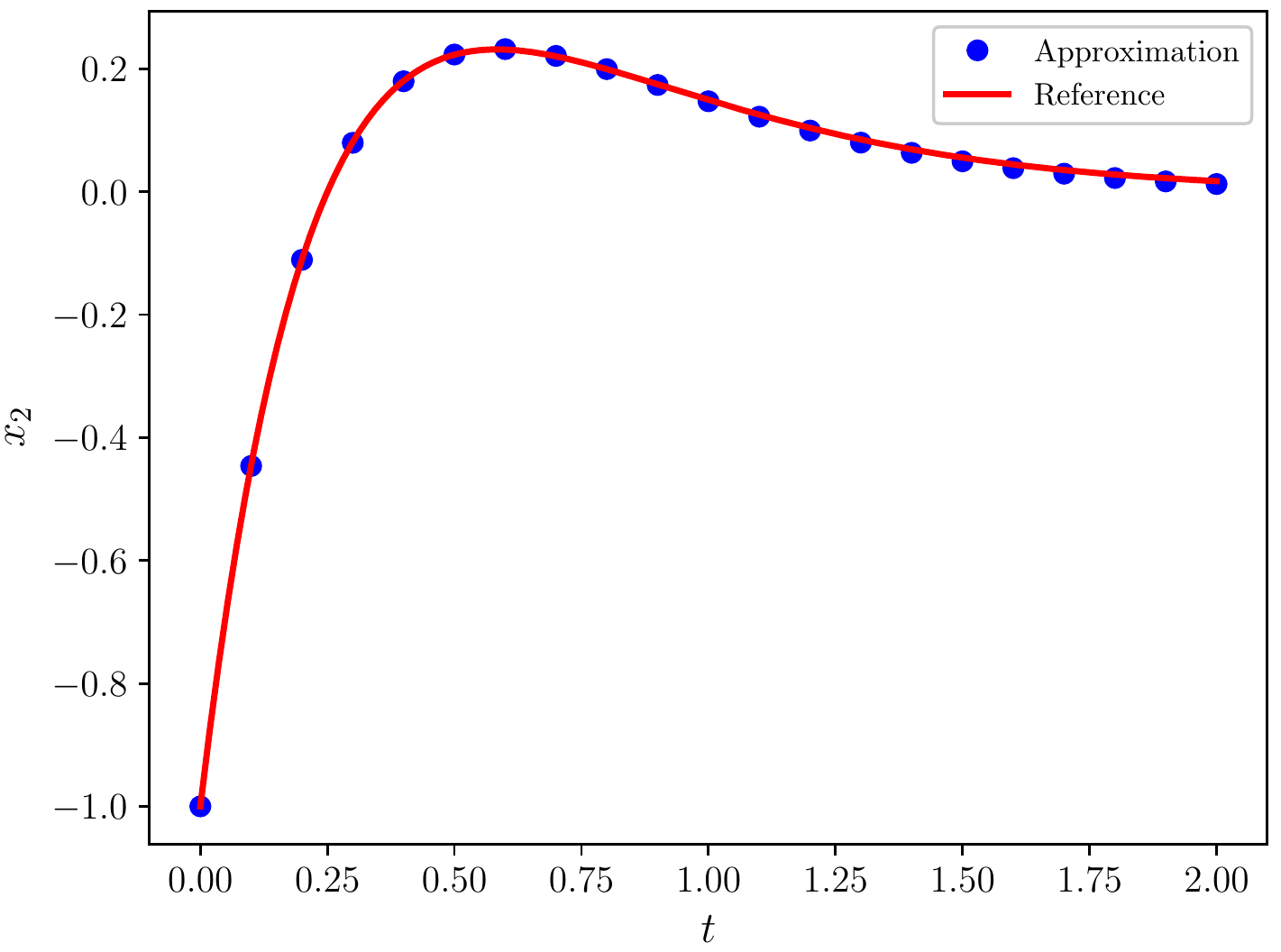}
		\caption{$x_2$, RT-ResNet}
	\end{subfigure}
	\begin{subfigure}[b]{0.32\textwidth}
	\begin{center}
		\includegraphics[width=1.0\linewidth]{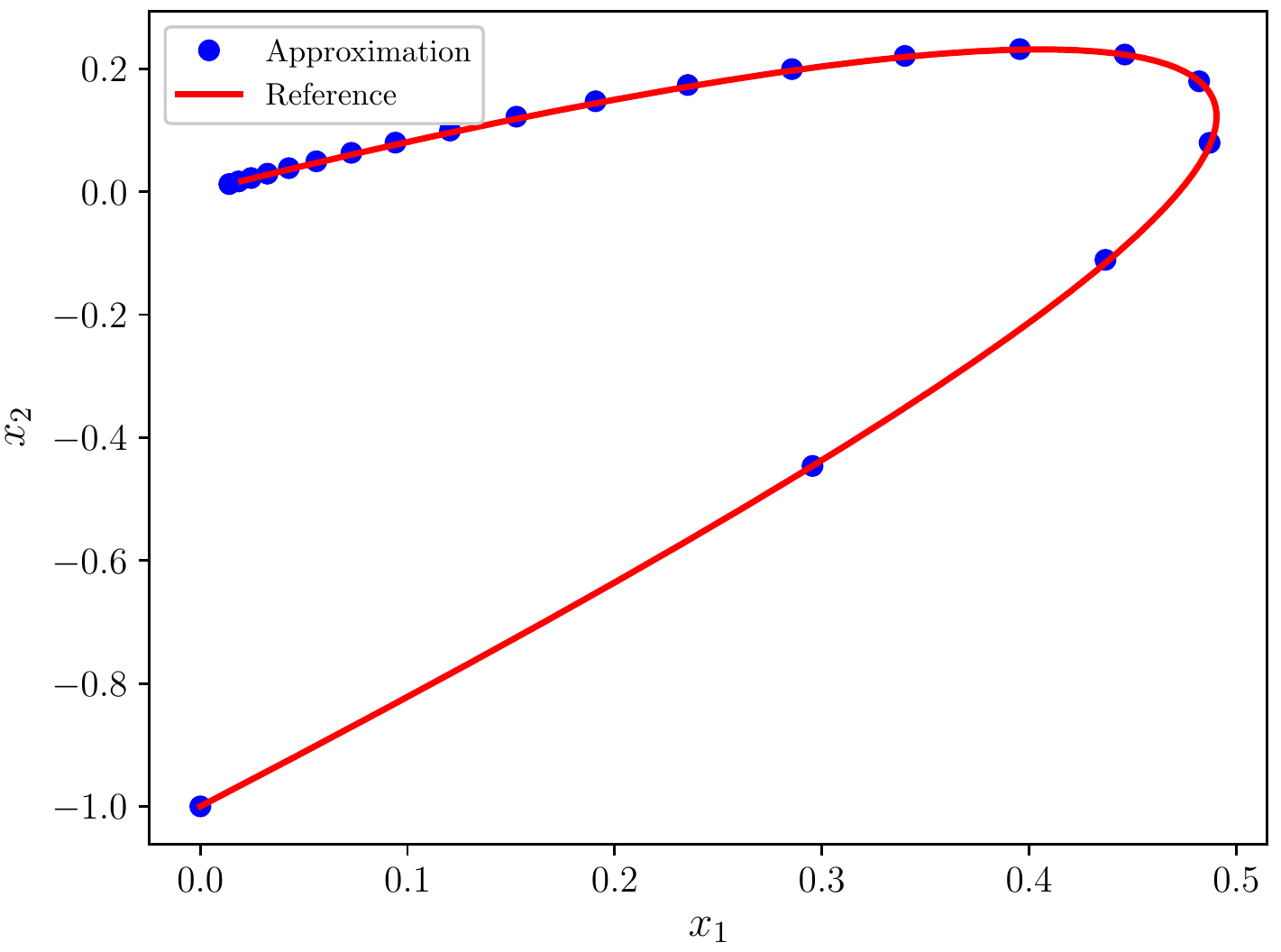}
		\caption{phase plot, RT-ResNet}
	\end{center}
	\end{subfigure}
	\begin{subfigure}[b]{0.32\textwidth}
	\begin{center}
		\includegraphics[width=1.0\linewidth]{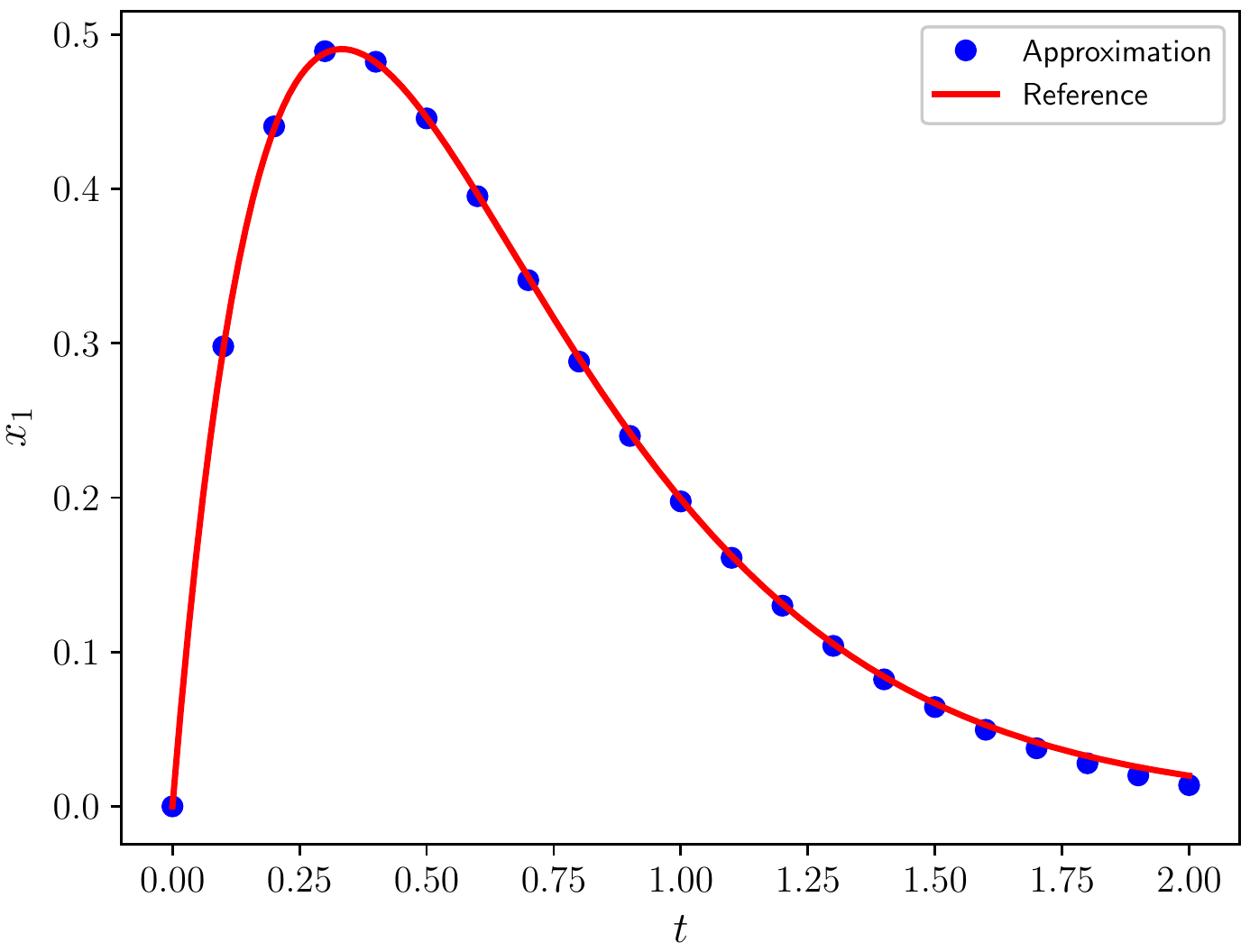}
		\caption{$x_1$, RS-ResNet}
	\end{center}
	\end{subfigure}
	\begin{subfigure}[b]{0.32\textwidth}
		\centering
		\includegraphics[width=1.0\linewidth]{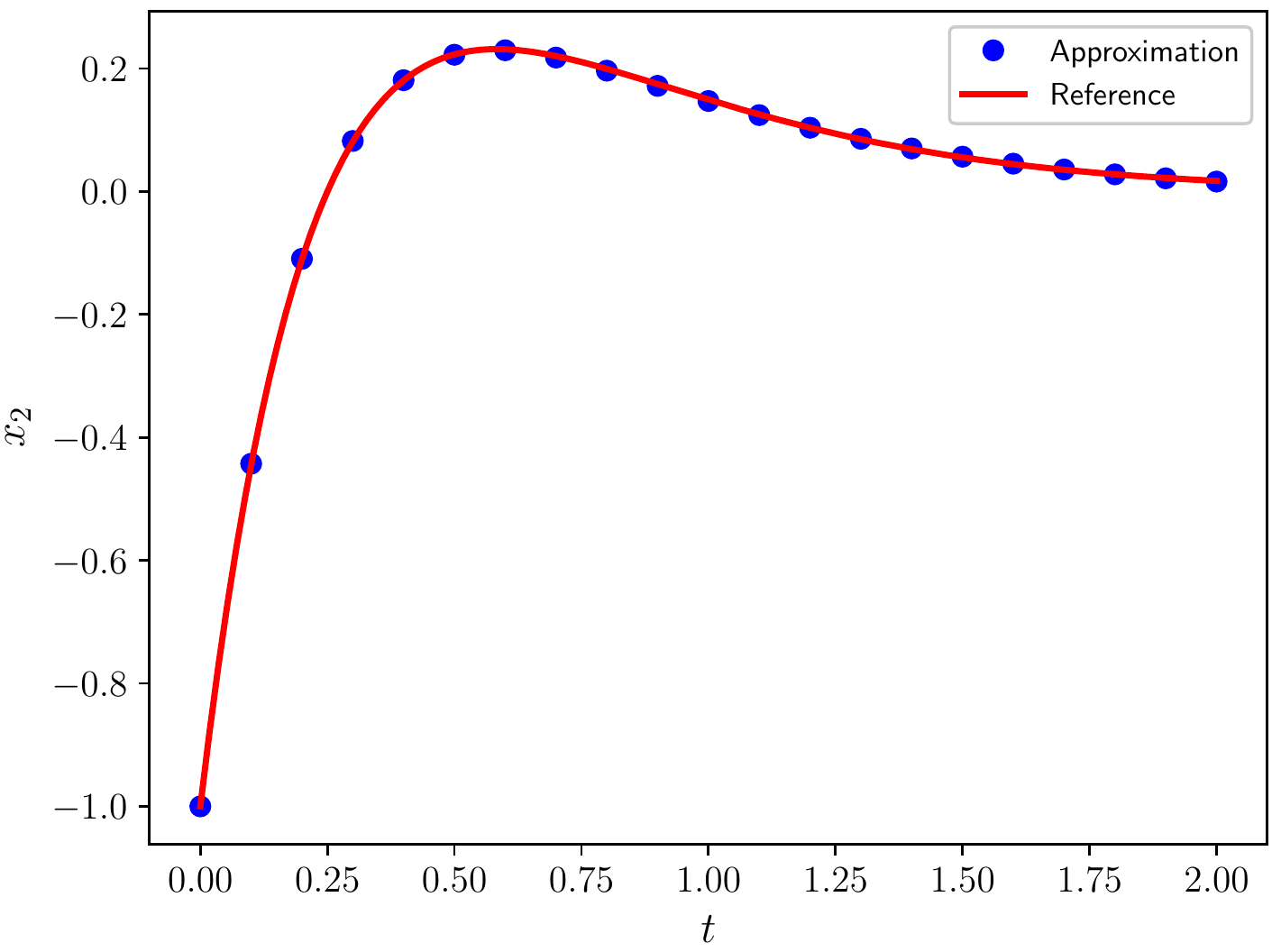}
		\caption{$x_2$, RS-ResNet}
	\end{subfigure}
	\begin{subfigure}[b]{0.32\textwidth}
	\begin{center}
		\includegraphics[width=1.0\linewidth]{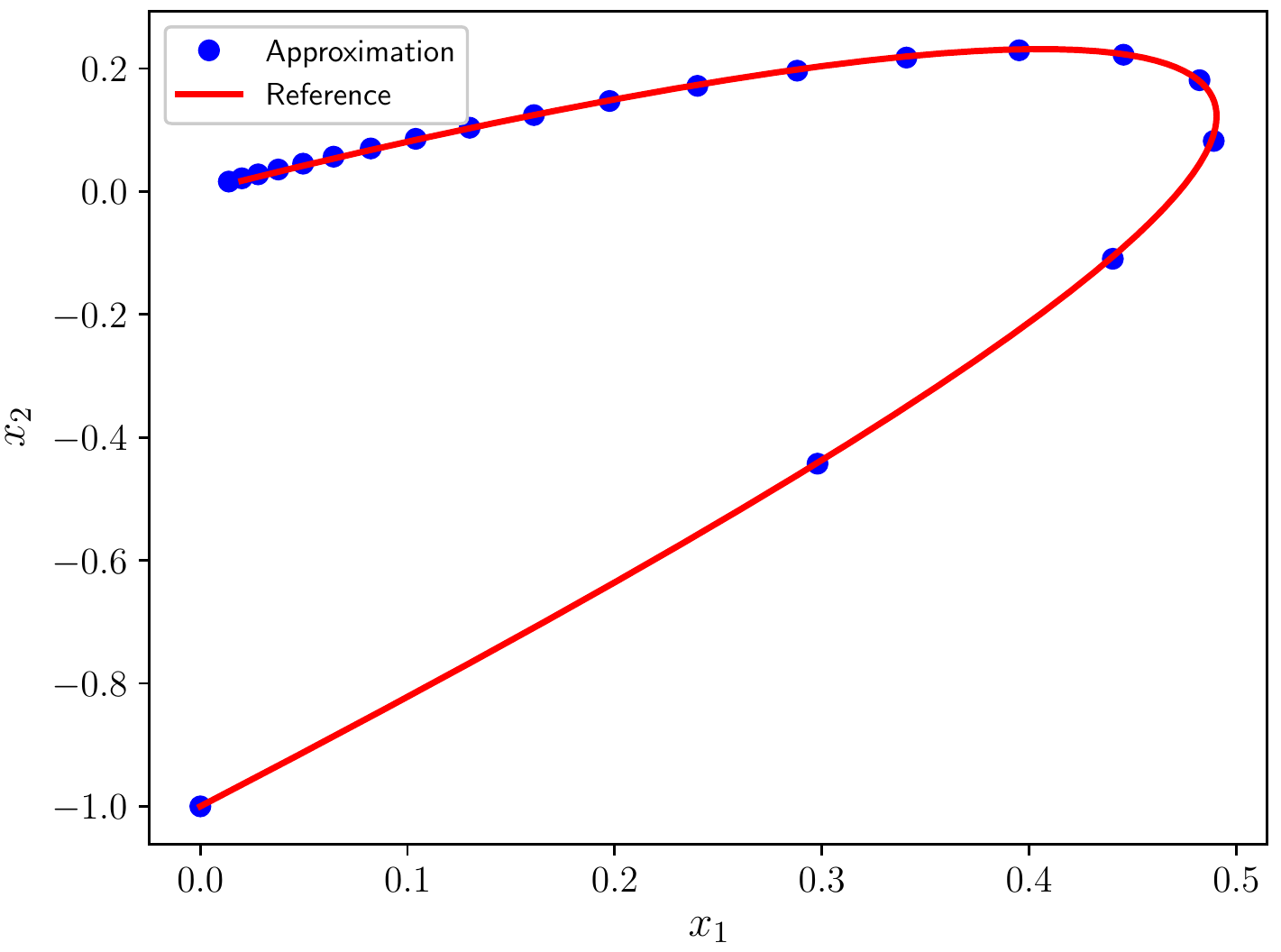}
		\caption{phase plot, RS-ResNet}
	\end{center}
	\end{subfigure}
	\caption{Trajectory and phase plots for the Example 2 with $\x_0=(0, -1)$.
	Top row: one-step ResNet model;
Middle row: Multi-step RT-ResNet model;
Bottom row: Multi-step RS-ResNet model.}
	\label{fig:ex2}
\end{figure}

We now consider another linear ODE system:
\begin{equation}
	\label{eq:example2}
	\begin{cases}
		\dot{x}_1=x_1-4x_2,\\
		\dot{x}_2=4x_1-7x_2.
	\end{cases}
\end{equation}
%We approximate the $\Delta=0.1$ flow map with the ResNet, the RT-ResNet and the
%RS-ResNet, respectively. The structure of each model is the same as in Example
%1. After learning the $\Delta$-lag flow
%map, we compose the DNN with itself for $20$ times to generate the
%approximation for the ODE solution. 
The numerical results for the three trained network models are
presented in Figure \ref{fig:ex2}.  Again, we show the prediction
results of the trained models for up to $t=2$. While all predictions
agree well with the reference solution, one can visually see that the
RS-ResNet model is more accurate than the RT-ResNet model, which in
turn is more accurate than the one-step ResNet model.
This is expected, as the multi-step methods should be more accurate
than the one-step method (ResNet), and the RS-ResNet should be even
more accurate due to its adaptive nature.
% We see that
%the approximation generated by the single-step ResNet, due to its simpler
%structure, shows slight deviations from the reference solution as the solution
%marches in time. For the three-step recurrent and recursive ResNets,
%since they are able to capture the more detailed structures in the single
%$\Delta$ step, the approximations stay much closer to the reference solution.
%Between the RS-ResNet and the RT-ResNet, the RS-ResNet gives better
%approximation, especially for the $x_1$ component. This is reasonable since the
%RS-ResNet has an ``adaptive'' nature and hence a better approximation ability. 

\subsection{Nonlinear ODEs}

We now consider two nonlinear problems. The first one is the well
known damped pendulum problem, and the second one is an nonlinear
differential-algebraic equation (DAE) for modelling a generic toggle
(\cite{gardner2000construction}).
In both examples, our one-step ResNet model has 2 hidden layers, each
of which has 40 neurons. Our multi-step RT-ResNet and RS-ResNet models
both have 3 of the same ResNet blocks ($K=3$). Again, our training
data are collected over $\Delta=0.1$ time lag. We produce predictions
of the trained model over time for up to $t=20$ and compare the
results against the reference solutions.

\subsubsection*{Example 3: Damped pendulum}

\begin{figure}[!htb]
	\centering
	\begin{subfigure}[b]{0.48\textwidth}
	\begin{center}
		\includegraphics[width=1.0\linewidth]{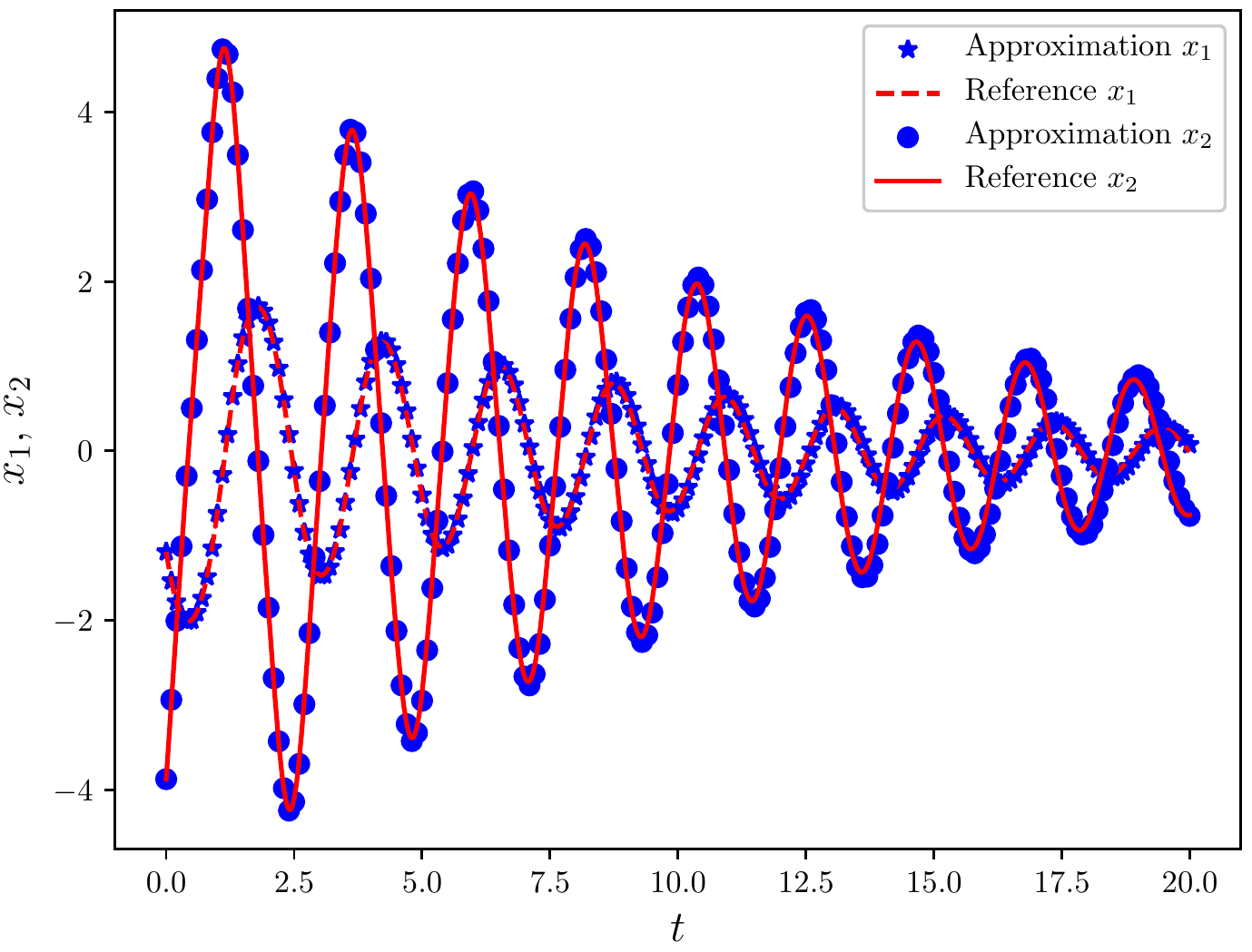}
		\caption{$x_1$ and $x_2$, ResNet}
	\end{center}
	\end{subfigure}
	\begin{subfigure}[b]{0.48\textwidth}
	\begin{center}
		\includegraphics[width=1.0\linewidth]{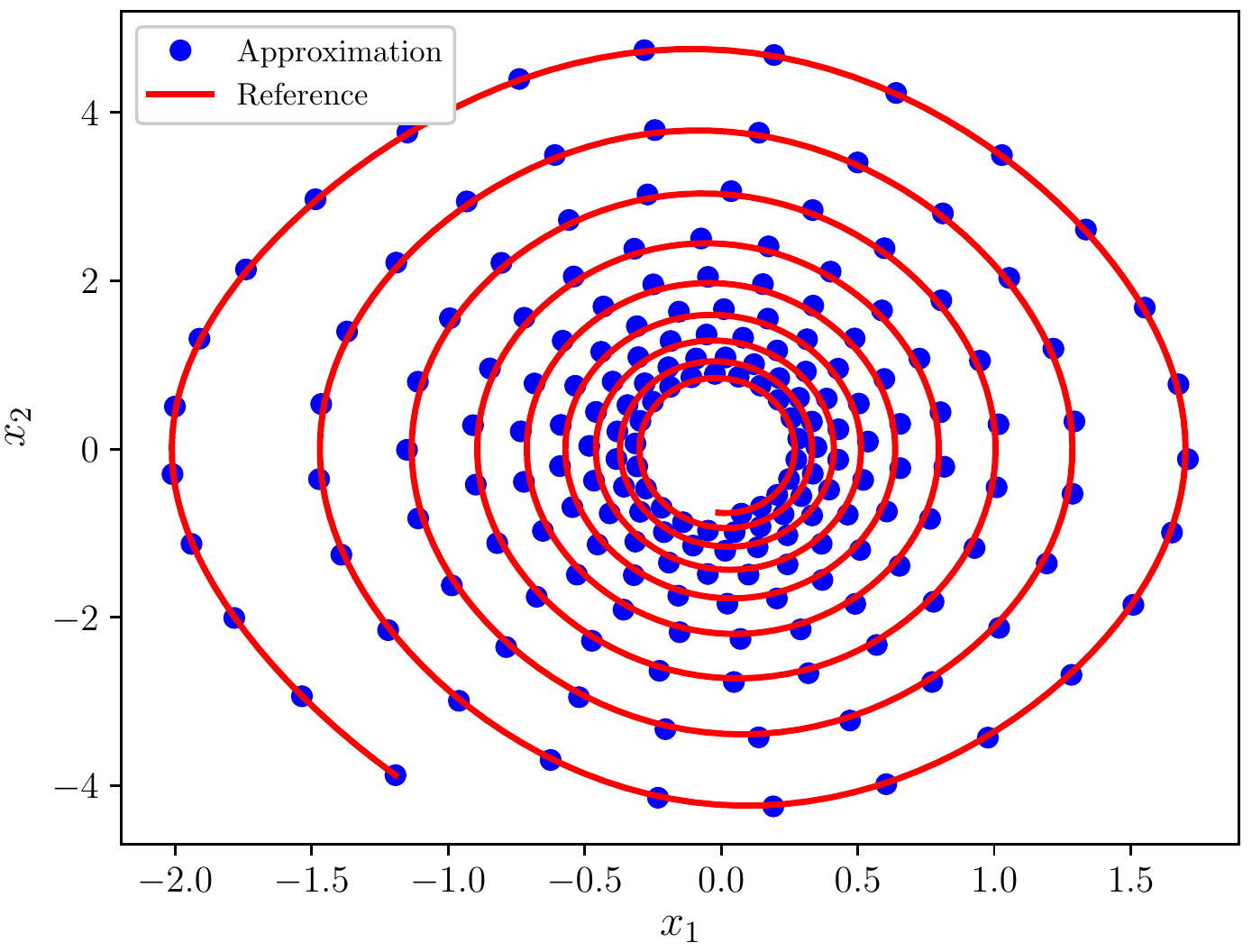}
		\caption{phase plot, ResNet}
	\end{center}
	\end{subfigure}
	\begin{subfigure}[b]{0.48\textwidth}
	\begin{center}
		\includegraphics[width=1.0\linewidth]{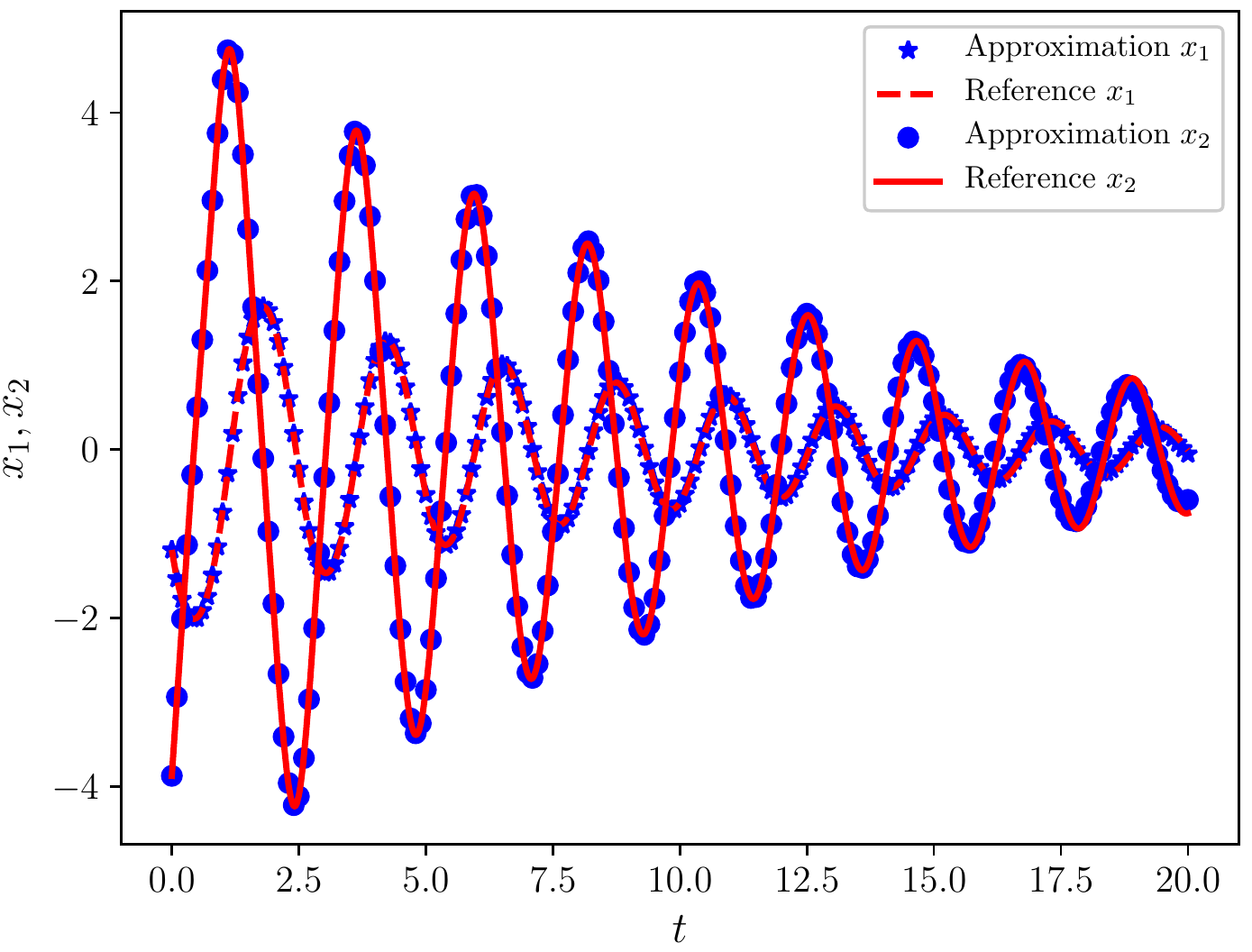}
		\caption{$x_1$ and $x_2$, RT-ResNet}
	\end{center}
	\end{subfigure}
	\begin{subfigure}[b]{0.48\textwidth}
	\begin{center}
		\includegraphics[width=1.0\linewidth]{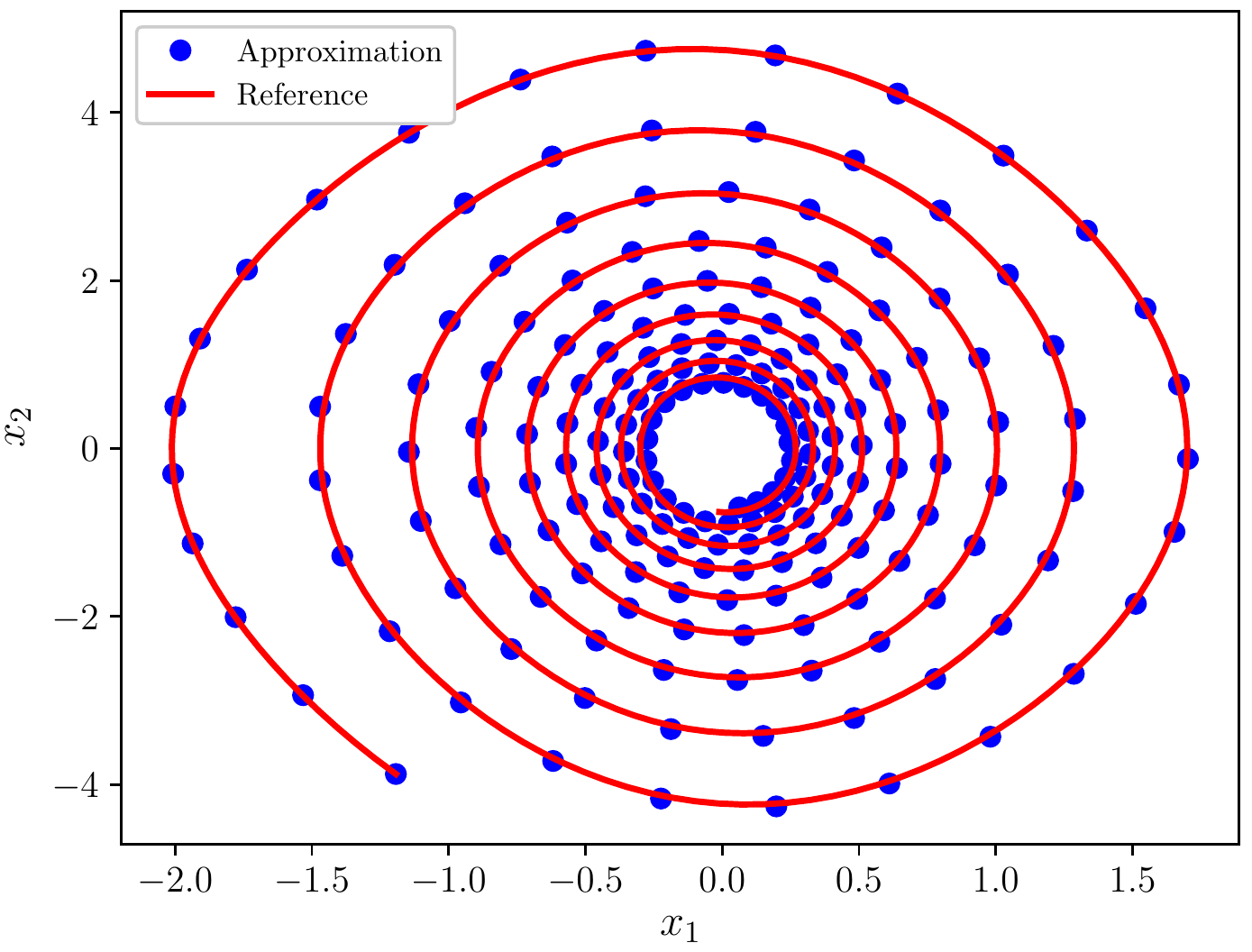}
		\caption{phase plot, RT-ResNet}
	\end{center}
	\end{subfigure}
	\begin{subfigure}[b]{0.48\textwidth}
	\begin{center}
		\includegraphics[width=1.0\linewidth]{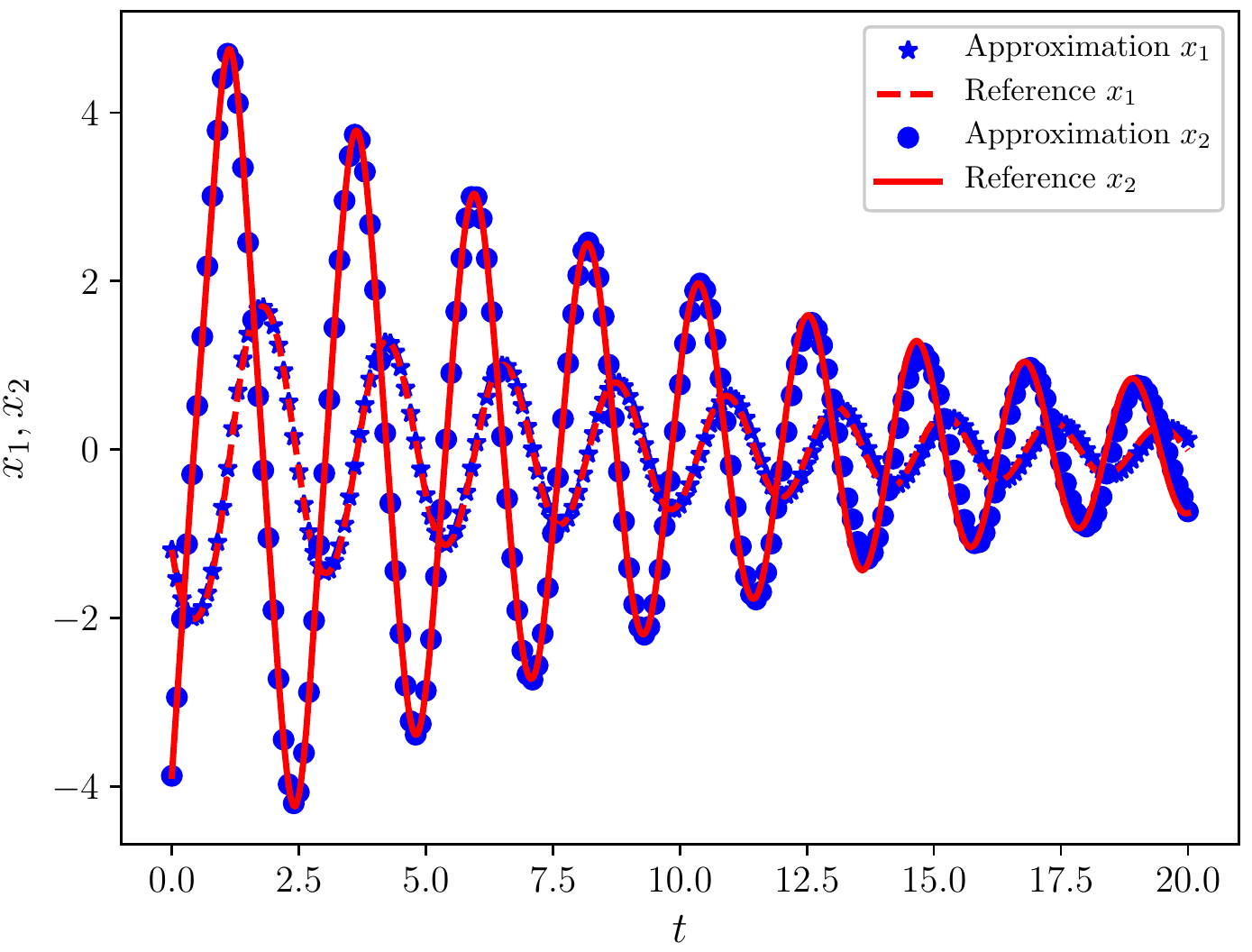}
		\caption{$x_1$ and $x_2$, RS-ResNet}
	\end{center}
	\end{subfigure}
	\begin{subfigure}[b]{0.48\textwidth}
	\begin{center}
		\includegraphics[width=1.0\linewidth]{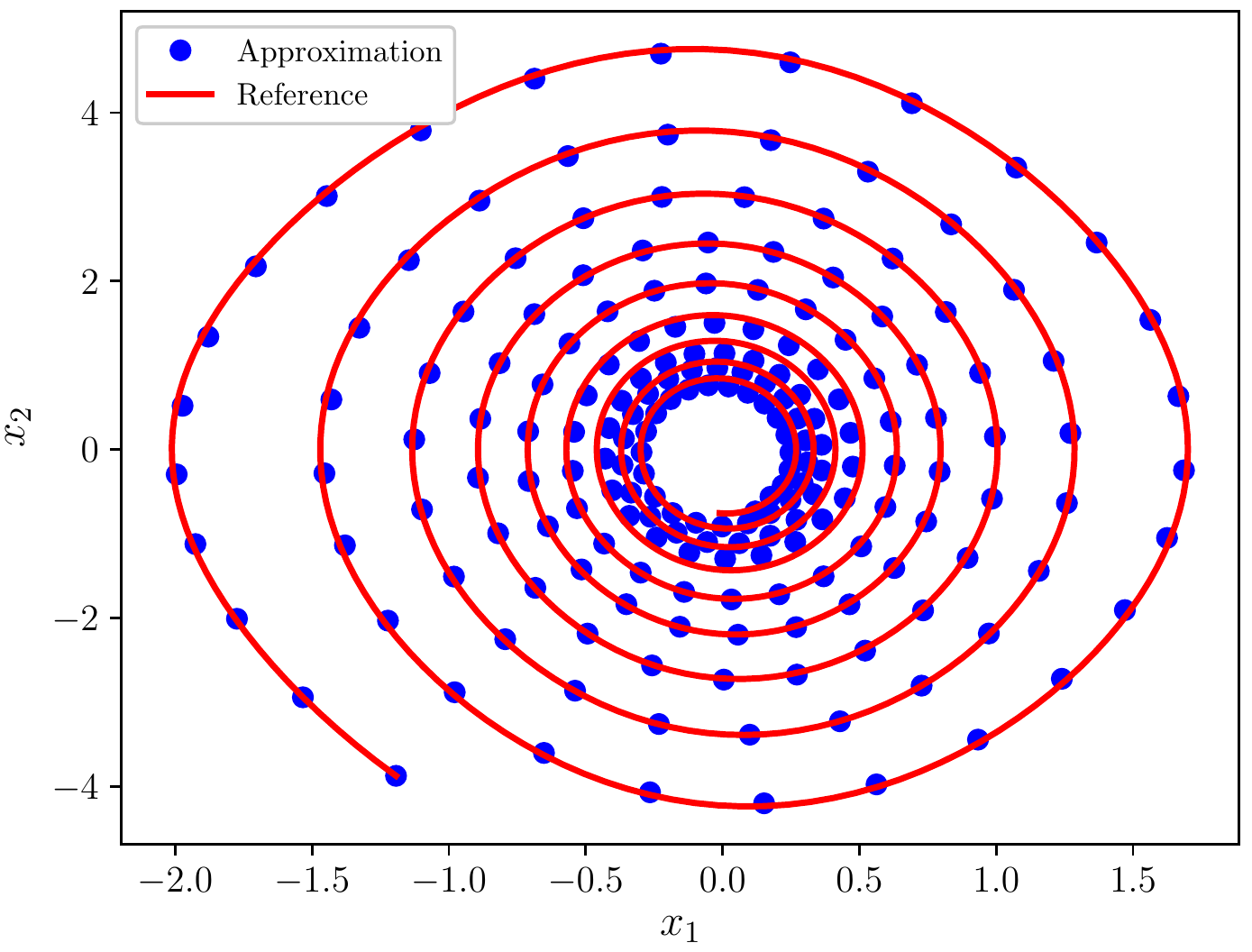}
		\caption{phase plot, RS-ResNet}
	\end{center}
	\end{subfigure}
	\caption{Trajectory and phase plots for the Example 3 with $\x_0=(-1.193, -3.876)$.
	Top row: one-step ResNet model;
Middle row: Multi-step RT-ResNet model;
Bottom row: Multi-step RS-ResNet model.}
	\label{fig:ex3}
\end{figure}

The first nonlinear example we are considering is the following damped
pendulum problem,
\begin{equation*}
	\begin{cases}
		\dot{x}_1=x_2,\\
		\dot{x}_2=-\alpha x_2-\beta \sin x_1.
	\end{cases}
\end{equation*}
where $\alpha=8.91$ and $\beta=0.2$. The computational domain is $D=[-\pi,
\pi]\times [-2\pi, 2\pi]$. 
In Figure \ref{fig:ex3}, we present the prediction results by the
three network models, starting from the initial condition $\x_0=(-1.193, -3.876)$ and
for time up to $t=20$.
We observe excellent agreements between the network models and the
reference solution.
\subsubsection*{Example 4: Genetic toggle switch}

We now considere a system of nonlinear differential-algebraic
equations (DAE), which are used to model a genetic toggle
switch in \textit{Escherichia coli}
(\cite{gardner2000construction}). It is composed of two repressors and
two 
constitutive promoters, where each promoter is inhibited by the represssor that
is transcribed by the opposing promoter. Details of experimental measurement
can be found in \cite{chartrand2011numerical}. This system of
equations are as follows,
\begin{equation*}
	\begin{cases}
		\dot{x}_1= \frac{\alpha_1}{1+x_2^\beta}-x_1,\\
		\dot{x}_2= \frac{\alpha_2}{1+z^\gamma}-x_2,\\
		z=\frac{x_1}{(1+\text{[IPTG]}/K)^\eta}.
	\end{cases}
\end{equation*}
In this system, the components $x_1$ and $x_2$ denote the concentration of the
two repressors. The parameters $\alpha_1$ and $\alpha_2$ are the effective
rates of the synthesis of the repressors; $\beta$ and $\gamma$ represent
cooperativity of repression of the two promoters, respectively; [IPTG] is
the concentration of IPTG, the chemical compound that induces the switch; and
$K$ is the dissociation constant of IPTG. 

In the following numerical experiment, we take $\alpha_1=156.25$,
$\alpha_2=15.6$, $\gamma=1$, $\beta=2.5$, $K=2.9618\times 10^{-5}$ and 
$\text{[IPDG]}=10^{-5}$. We consider the computational domain $D=[0,
20]^2$.  

\begin{figure}[!htb]
	\centering
	\begin{subfigure}[b]{0.48\textwidth}
	\begin{center}
		\includegraphics[width=1.0\linewidth]{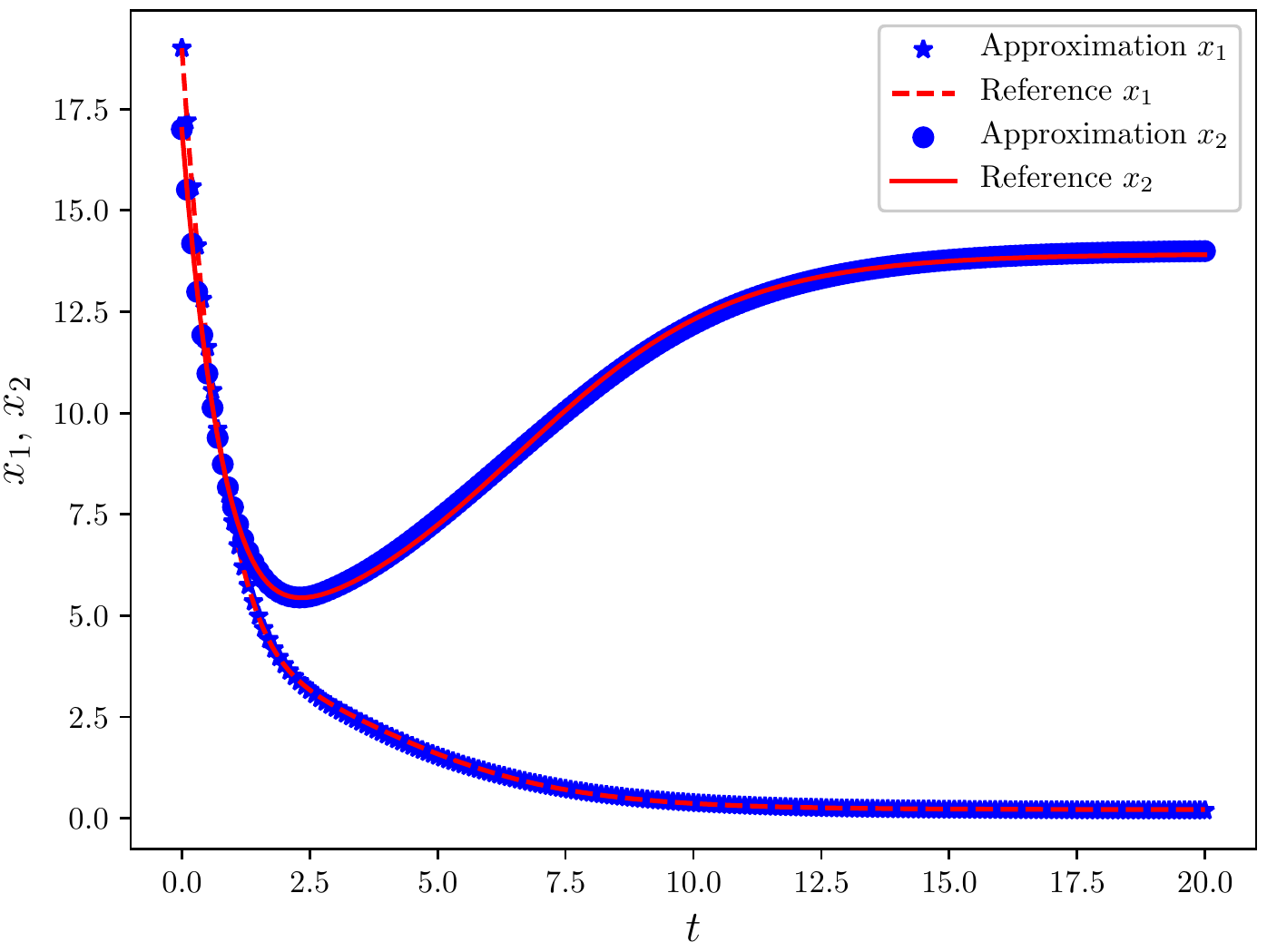}
		\caption{$x_1$ and $x_2$, ResNet}
	\end{center}
	\end{subfigure}
	\begin{subfigure}[b]{0.48\textwidth}
	\begin{center}
		\includegraphics[width=1.0\linewidth]{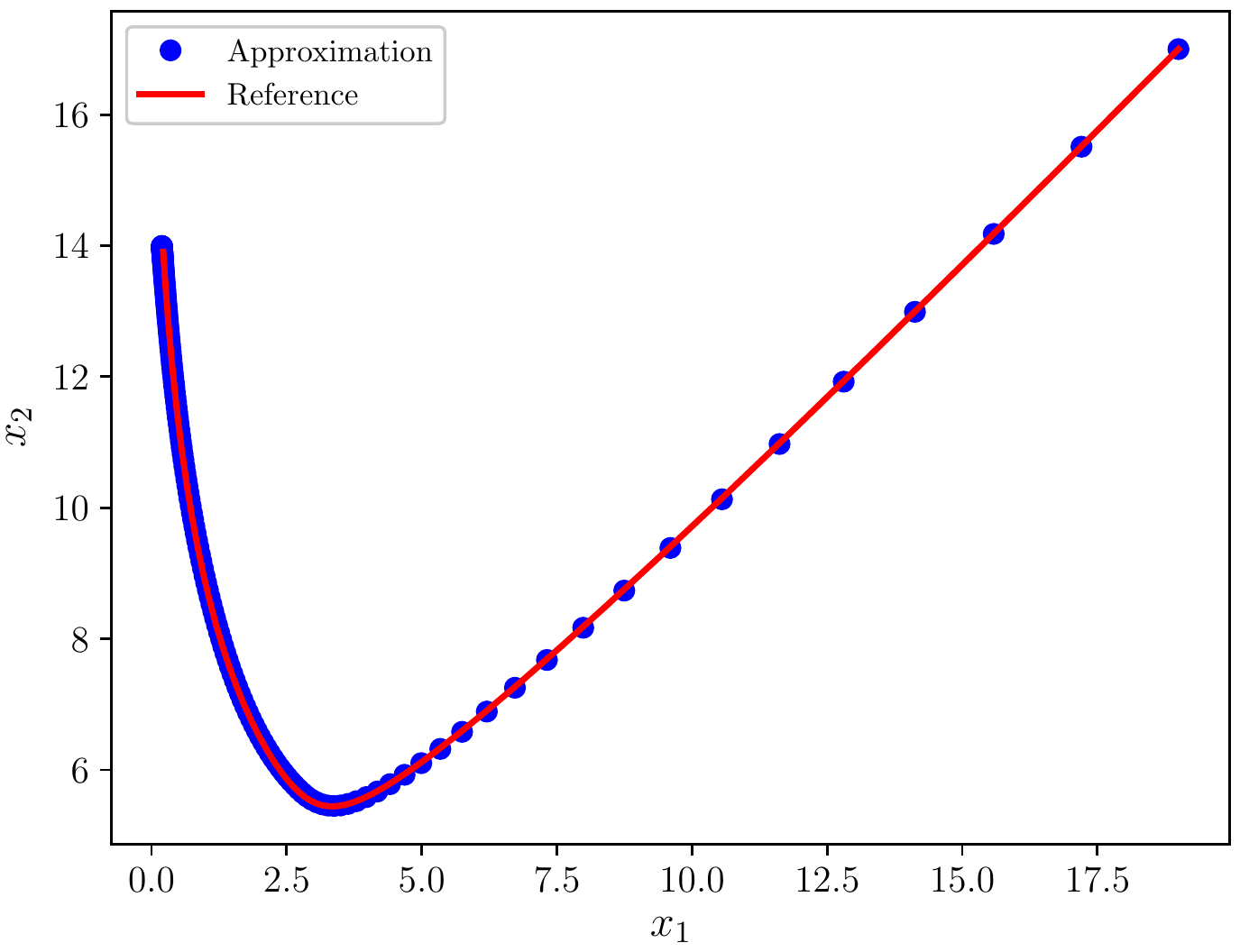}
		\caption{phase plot, ResNet}
	\end{center}
	\end{subfigure}
	\begin{subfigure}[b]{0.48\textwidth}
	\begin{center}
		\includegraphics[width=1.0\linewidth]{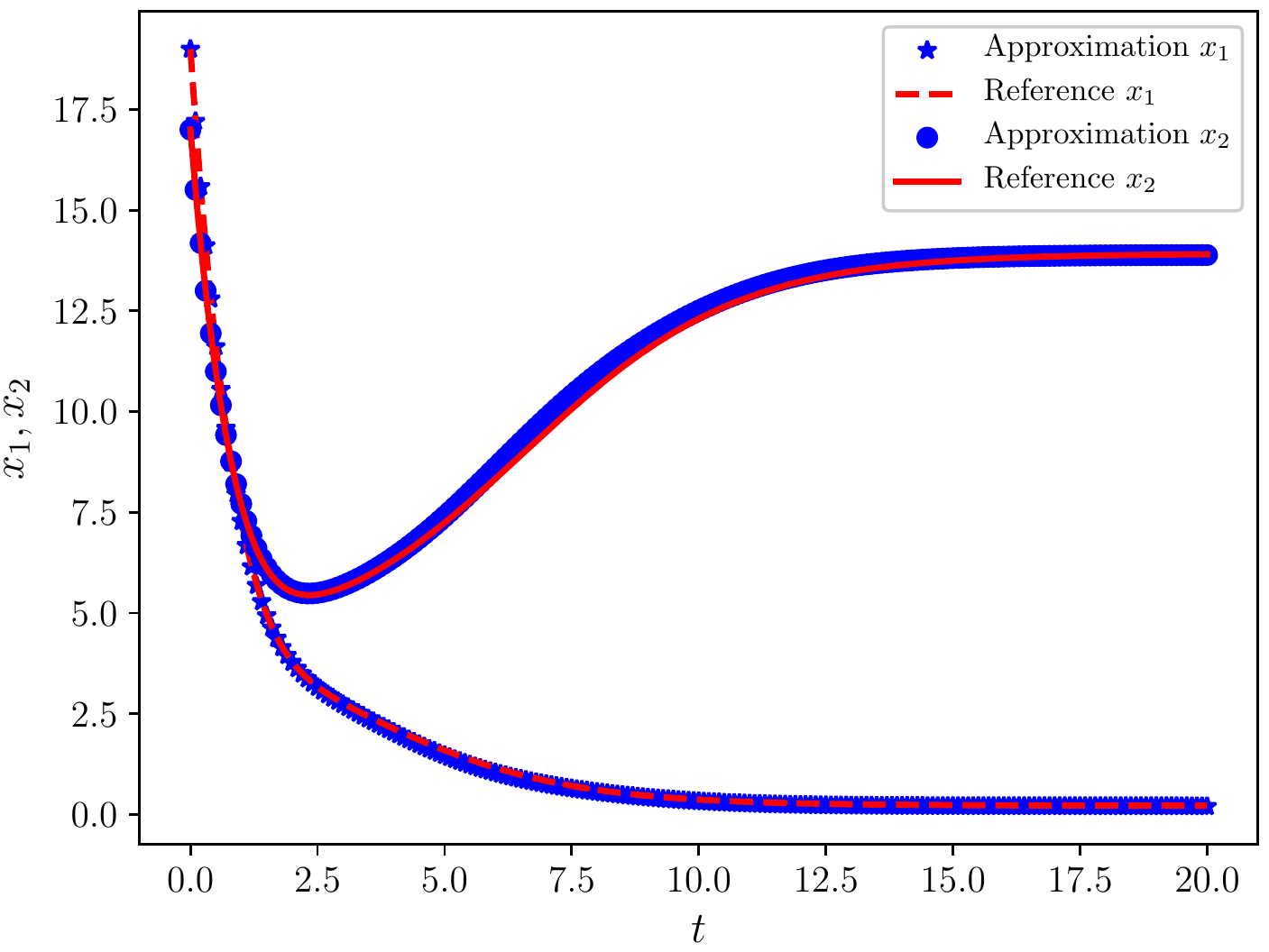}
		\caption{$x_1$ and $x_2$, RT-ResNet}
	\end{center}
	\end{subfigure}
	\begin{subfigure}[b]{0.48\textwidth}
	\begin{center}
		\includegraphics[width=1.0\linewidth]{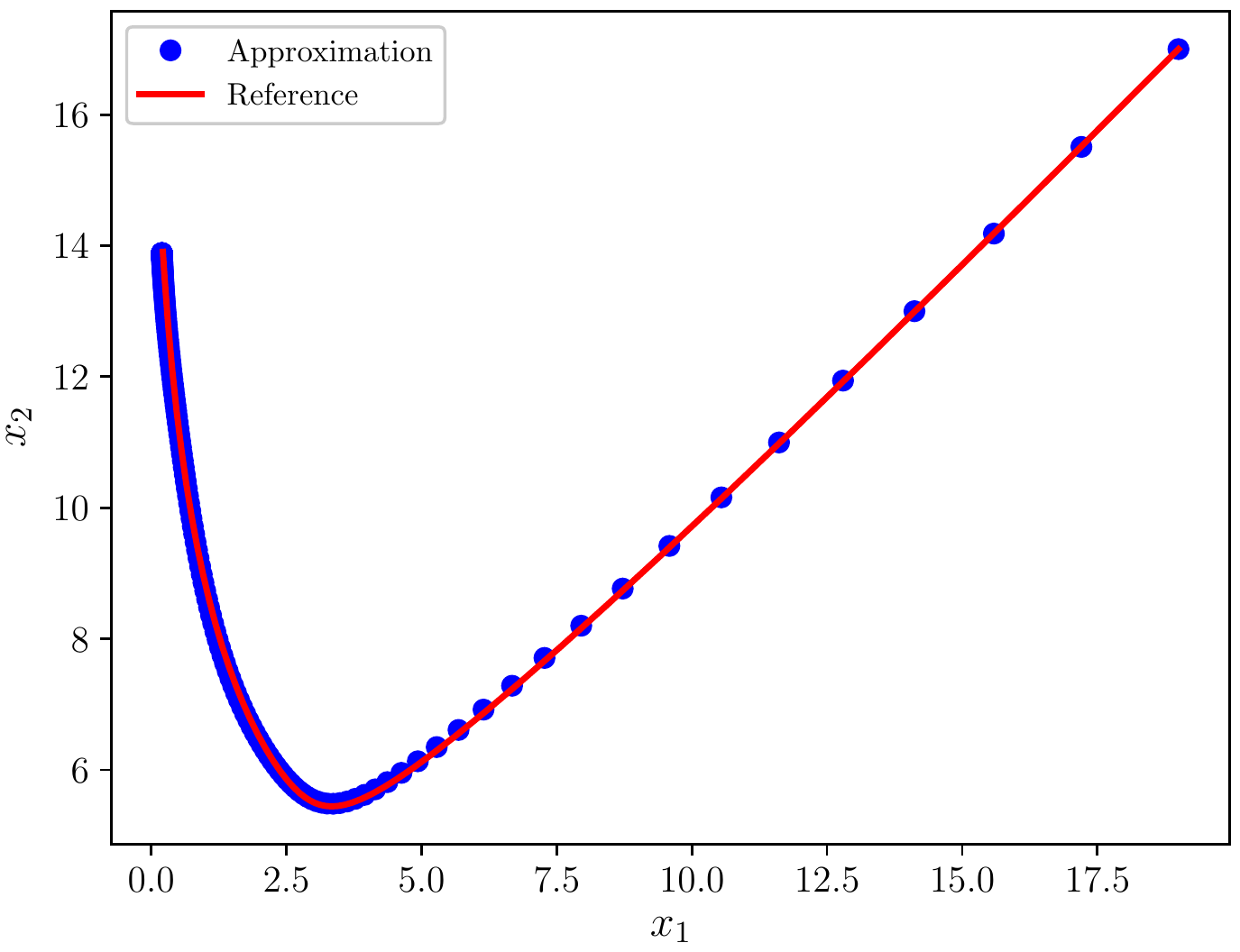}
		\caption{phase plot, RT-ResNet}
	\end{center}
	\end{subfigure}
	\begin{subfigure}[b]{0.48\textwidth}
	\begin{center}
		\includegraphics[width=1.0\linewidth]{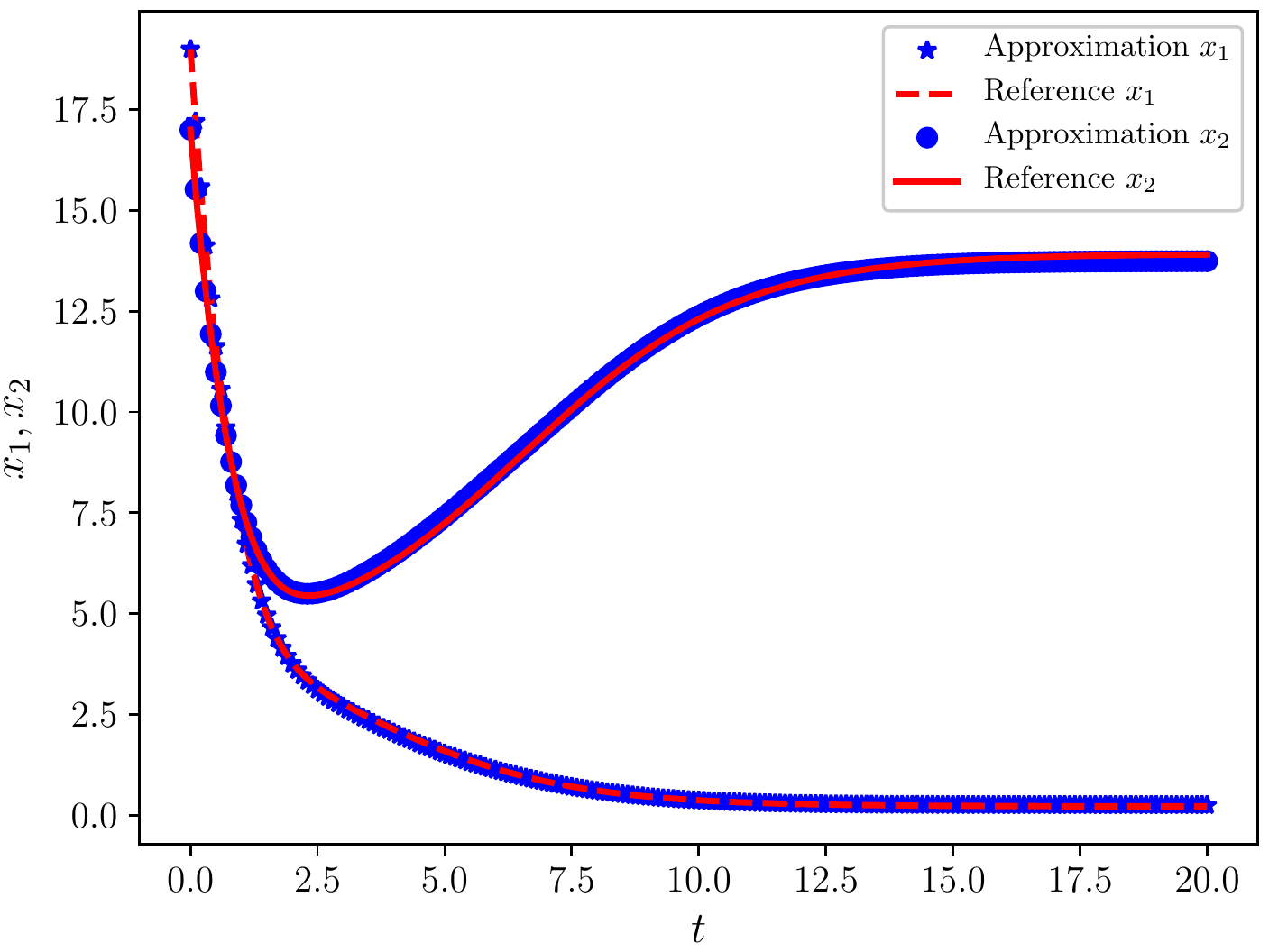}
		\caption{$x_1$ and $x_2$, RS-ResNet}
	\end{center}
	\end{subfigure}
	\begin{subfigure}[b]{0.48\textwidth}
	\begin{center}
		\includegraphics[width=1.0\linewidth]{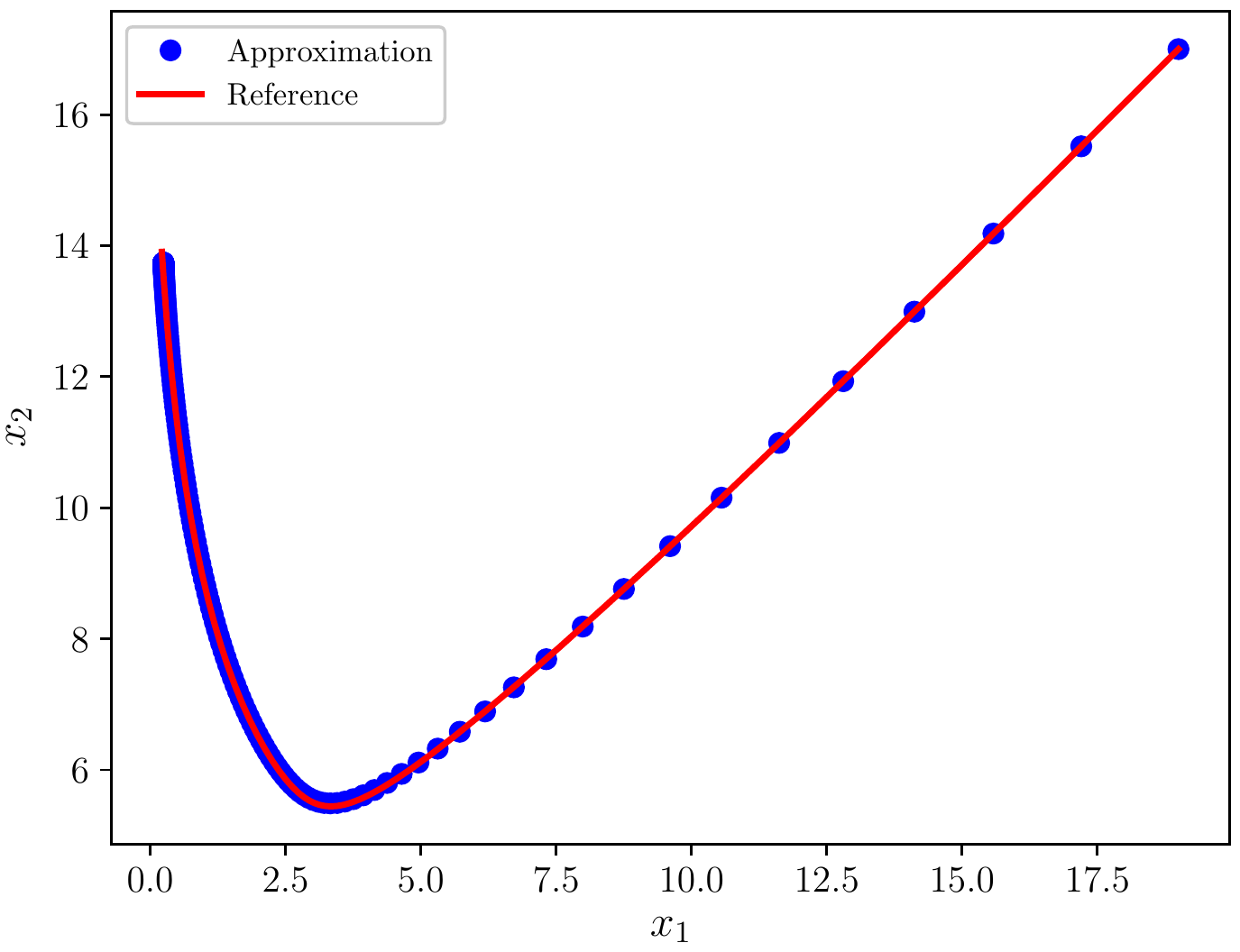}
		\caption{phase plot, RS-ResNet}
	\end{center}
	\end{subfigure}

	\caption{Trajectory and the phase plots for the Example 4 with $\x_0=(19, 17)$.
	Top row: one-step ResNet model;
Middle row: Multi-step RT-ResNet model;
Bottom row: Multi-step RS-ResNet model.}
	\label{fig:ex4}
\end{figure}

In Figure \ref{fig:ex4} we present the prediction results generated by the
ResNet, the RT-ResNet and the RS-ResNet, for time up to $t=20$. 
The initial condition is
$\x_0=(19, 17)$. Again, all these three models produce accurate
approximations, even for such a long-time simulation. 

\section{Conclusion} \label{sec:conclusions}

We presented several deep neural network (DNN) structures for
approximating unknown dynamical systems using trajectory data. The DNN
structures are based on residual network (ResNet), which is a one-step
method exact time integrator.
Two multi-step variations were presented. One is recurrent ResNet
(RT-ResNet) and
the other one is recursive ResNet (RS-ResNet). Upon successful
training, the methods produce discrete dynamical systems that
approximate the underlying unknown governing equations. All methods
are based on integral form of the underlying system. Consequently,
their
constructions do not require time derivatives of the trajectory data
and can work with coarsely distributed data as well.
We presented
the construction details of the methods, their theoretical
justifications, and used several examples to demonstrate the
effectiveness of the methods.

%\section*{Acknowledgment}
%This work is in part supported by AFOSR, DOE, NNSA, NSF.

\bibliographystyle{plain}
\bibliography{neural,LearningEqs}

\end{document}